\documentclass[11pt]{article}

\usepackage[a4paper,margin=1in]{geometry}
\usepackage{amsmath,amssymb,amsfonts,amsthm,mathtools}
\usepackage{enumitem,mathrsfs}
\usepackage{needspace}
\usepackage{xcolor}
\usepackage[T1]{fontenc}
\usepackage{dsfont}

\usepackage{mathtools}
\mathtoolsset{showonlyrefs}

\usepackage[sort,nocompress]{cite}

\usepackage{hyperref}

\hypersetup{
    colorlinks=true,
    linkcolor=blue,
    citecolor=blue,
    urlcolor=blue
}

\newcommand{\RR}{\mathbb R}

\newcommand{\ZZ}{\mathbb Z}
\newcommand{\NN}{\mathbb N}

\newcommand{\Prob}{\mathbb P}

\newcommand{\Tr}{\operatorname{Tr}}

\newcommand{\E}{\operatorname{E}}
\newcommand{\Var}{\operatorname{Var}}

\newcommand{\rank}{\operatorname{Rank}}

\newcommand{\Sone}{\mathcal S_1}
\newcommand{\Stwo}{\mathcal S_2}
\newcommand{\mbbH}{\mathbb{H}}

\newcommand{\supp}{\operatorname{supp}}
\allowdisplaybreaks[2]
\newcommand{\dd}{\,d}

\usepackage[ruled]{algorithm2e}
\usepackage{lscape}
\usepackage{rotating}
\numberwithin{equation}{section}
\newtheorem{theorem}{Theorem}[section]
\newtheorem{proposition}{Proposition}[section]
\newtheorem{lemma}{Lemma}[section]
\newtheorem{corollary}{Corollary}[section]

\theoremstyle{definition}
\newtheorem{remark}{Remark}[section]

\newtheorem{assumption}{Assumption}[section]

\newcommand{\mockalph}[1]{}

\title{A Dobrushin--Major--Taqqu Theorem for Hilbert Space-Valued Random
Variables}
\author{
Marie-Christine D\"uker 
\\ Technical University of Munich
\and
María Dolores Ruiz-Medina
\\ University of Granada
\and
Pavlos Zoubouloglou      
\\ University of Münster
}
\date{\today}

\begin{document}

\maketitle

\begin{abstract}
Let $\{X_t\}_{ t \in \ZZ}$ be a stationary Gaussian process taking values in a separable Hilbert space $\mathbb H_1$, with long-range dependence defined in the spectral domain through a memory operator $A:\mathbb H_1 \to \mathbb H_1$. For square-integrable, nonlinear operators $G:\mathbb H_1\to\mathbb H_2$, taking values in another separable Hilbert space $\mathbb H_2$, we  establish a non-central limit theorem for the centered partial-sum process of $\{G[X_t]\}_{t \in \ZZ}$, whenever $G[X_t]$ inherits the long-range dependence behavior from its Gaussian input process. 
We identify an operator normalization that accounts for both the underlying memory structure, characterized by the memory operator $A$, and the nonlinear transformation $G$. The limiting Hilbert space-valued process is, in general, non-Gaussian and can be expressed as multiple Wiener integrals involving only finitely many Wiener chaos orders along fixed directions in $\mathbb{H}_2$.
The proof introduces a Hilbert space-valued generalization of the reduction principle, identifying the dominant Hermite components and showing that the remainder is negligible under our operator normalization.
We employ our results to prove non-central limit theorems for nonlinear transformations of functional FARIMA processes and sojourn functionals of temporally dependent 
spherical Gaussian fields.
\end{abstract}

\section{Introduction}
Let \(\{X_t\}_{t\in\ZZ}\) be a centered stationary Gaussian sequence with values in a separable Hilbert space \(\mathbb H_1\) and \(G:\mathbb H_1\to\mathbb H_2\) a measurable map, with \(\mathbb H_2\) another separable Hilbert space. Define the partial-sum process
\begin{equation}
    \label{eq:partial-sums}
    S_n(s)
    =
    T_n
    \sum_{t=1}^{\lfloor ns\rfloor}
    \left(G[X_t]-\E G[X_t]\right),
    \qquad s\in[0,1],
\end{equation}
where \(\{T_n\}_{n \geq 1}\) is a suitable sequence of scaling operators. This work aims to understand the weak limit of \(\{S_n\}_{n \geq 1}\) in \eqref{eq:partial-sums} when the process \(\{G[X_t]\}_{t \in \ZZ} \) exhibits long-range dependence.

The complementary regime, in which \(\{G[X_t]\}_{ t \in \ZZ}\) exhibits short-range dependence, is customarily referred to as the Breuer--Major setting, in recognition of the seminal work \cite{Breuer83} in the case \(\mathbb H_1 = \mathbb H_2 = \RR\); see also \cite{GiraitisSurgailis1985}. In this case, the process \eqref{eq:partial-sums} falls within the Gaussian universality class, with the usual scaling \(T_n=n^{-1/2}I_{\mathbb H_2}\). More recently, a line of work based on Malliavin calculus and Stein's method has provided new proofs and extensions of such central limit theorems, e.g., \cite{NourdinPeccati2012,NourdinNualart2020,CampeseNourdin20}. Using the Malliavin--Stein framework in infinite-dimensional spaces \cite{DukerZoub25}, two of the present authors recently established the first extension of the Breuer--Major theorem to the present framework, formulated entirely in Hilbert spaces \cite{DukerZoub26}.

In the long-range dependent regime, the limit of \eqref{eq:partial-sums} is generally non-Gaussian and such weak convergence results are referred to as non-central limit theorems. In the scalar case $\mbbH_1 = \mbbH_2 = \mathbb{R}$, the corresponding non-central limit theorems are due to Dobrushin and Major \cite{Dobrushin} and Taqqu \cite{Taqqu,Taqqu2}, whose {work motivates} the title of the present article. In this case, the scaling $T_n$ depends, in general, on the dependence properties of $\{G(X_t)\}_{ t \in \ZZ}$, and the partial sums $\{S_n\}_{n \geq 1}$ fall under the universality classes of the so-called Hermite processes. The cases $\mbbH_1 = \mathbb{R}^d, \mbbH_2 = \mathbb{R}$ and $\mbbH_1 = \mathbb{R}, \mbbH_2 = \mathbb{R}^d$ were later studied in \cite{SanchezDeNaranjo1993,Arcones1994} and \cite{BaiTaqqu}, respectively. This topic continues to receive attention; see, for instance, \cite{Major2026Noncentral} for the case $\mbbH_1 = \mathbb{R}^d, \mbbH_2 = \mathbb{R}$, with a latent Gaussian random field $\{X_t\}_{t \in \mathbb{Z}^r}$ indexed by the lattice $\mathbb{Z}^r$. Finally, Theorem~2.7 of \cite{BetkenDehling2021} studies the case $\mbbH_1 = \mathbb{R}, \mbbH_2 = L^2(\mathbb{R},w(t)dt)$, where $w$ is a suitable weight function, and $T_n$ is a scalar multiple of the identity times a slowly varying function, with the scalar normalization determined by the leading term. Their result can be recovered as a special case of ours, modulo the presence of the slowly varying function, giving rise to a substantially simpler limiting process.

We establish a non-central limit theorem allowing both the input and output spaces to be general, possibly infinite-dimensional Hilbert spaces. Under spectral assumptions on $\{X_t\}_{ t \in \ZZ}$ (cf. Assumption \ref{ass:trace-admissibility}) and mild assumptions on the map $G$ (cf. Assumption \ref{ass:bigAssump2}), we determine the limit of \eqref{eq:partial-sums}, with \(\{T_n\}_{n \geq 1}\) a suitable sequence of scaling operators. In this setting, different projections in $\mbbH_2$ may require different normalizations, and different Wiener chaos orders may contribute to the limit, even within a single direction.
We identify the operator normalization \(T_n=n^{-D}\) through the effective memory exponents of the active Hermite terms and represent the limiting process by a series of Hilbert-valued multiple Wiener integrals.

Our motivation for this problem stems from the asymptotic analysis of sojourn measures of Gaussian random fields exhibiting long-range dependence. To be more precise, consider a Gaussian random field \(Z=\{Z(x,t):x\in\mathcal D,\ t\in\mathcal T\}\), where \(\mathcal D\) is a spatial domain such as \([0,T]^d\) or the sphere \(\mathbb S^2\) \cite{Marinucciv2}, \(\mathcal T\) is a time index set such as \(\mathbb R\) or \(\mathbb Z\), and \(\{Z(\cdot,t)\}_{t\in\mathcal T}\) is stationary and exhibits long-range dependence. Viewing \(Z(\cdot,t)\) as an \(L^2(\mathcal D)\)-valued random element, set \(\mbbH_1=L^2(\mathcal D)\) and \(\mbbH_2=\mathbb R\). For a measurable domain \(D\subseteq\mathcal D\) of finite measure, the centered sojourn functional above \(u\) is, for $k\in\mathcal T$,
\begin{equation}
    \begin{aligned}
        \mathcal S_D^k(u)
        = G_u[Z(\cdot,k)]-\E G_u[Z(\cdot,0)], \quad
        G_u[f]
        = \int_D \mathds{1}_{\{f(x)\ge u\}}\,dx
         = \left|\{x\in D:f(x)\ge u\}\right|,
    \end{aligned}
\end{equation}
where \(|\cdot|\) denotes Lebesgue or surface measure, as appropriate. The associated asymptotic problem concerns the normalization and limiting distribution of centered sojourn measures as the observation domain expands. For a time-indexed spatial field, this also leads to the study of temporal sums of \(G_u[Z(\cdot,t)]-\E G_u[Z(\cdot,0)]\).

These questions have been investigated under both short- and long-range dependence; see \cite{Bulinski,Ivanov89,Leonenko99,LeonenkoOlenko14,Makogin20,Marinucciv2}. A central tool in this literature is the reduction principle, which identifies the chaos components governing the asymptotic fluctuations through Hermite or Laguerre expansions. Classical results \cite{Dobrushin,Taqqu,Taqqu2,Leonenkoetal17,Leonenkoetal17b} establish this approach for long-range dependent Gaussian and chi-squared subordinated fields. In the spatiotemporal setting, Leonenko and Ruiz--Medina \cite{LeonenkoRuiz-Medina2023,Leonenkoruizmedi25} apply reduction principles to obtain central limit theorems under increasing-domain asymptotics in time and/or space. Marinucci, Rossi and Vidotto \cite{MarinucciRV} derive central and non-central limit theorems for sojourn measures of stationary and spatially isotropic Gaussian fields on the sphere under increasing-time asymptotics, showing how temporal dependence and the excursion threshold determine the normalization and limiting law. Their spherical-harmonic formulation makes the role of the angular spectrum explicit. Extensions to chi-squared subordinated fields, using tensor products of Hermite polynomial bases, are developed in \cite{Caponera}. See Section \ref{subsec:spherical-excursions} for the application of our result to such a setting.

Although we are the first to treat the partial sums of subordinated Gaussian sequences \eqref{eq:partial-sums} at this level of generality, several works establish limit theorems for stationary Hilbert-valued sequences. The theory is more developed for linear series \cite{Bosq2000,Merlevede1996,MerlevedePeligradUtev1997}, and a few works also study limit theorems for long-range dependent linear series in Hilbert spaces \cite{CharaciejusRackauskas2014,Duker2018,DukerZoubouloglou2025}. In particular, to the best of our knowledge, the last of these works provides the first construction of a non-Gaussian Hermite process with values in an infinite-dimensional Hilbert space.

Our main result, Theorem~\ref{thm:operator-normalized-full}, establishes a non-central limit theorem under spectral assumptions, with the anticipated operator normalization \(n^{-D}\) on \(\mathbb H_2\). A key assumption is that a single orthonormal basis of \(\mbbH_1\), denoted by \(\{\phi_i\}_{i\in\NN}\), simultaneously diagonalizes all covariance and spectral density operators. Moreover, we assume a semiparametric characterization of the spectral densities in this basis, such that each coordinate exhibits long-range dependence, with uniform control of the associated memory parameters and slowly varying functions; see Assumption~\ref{ass:trace-admissibility}. 
Our assumptions on the functional time series \(\{X_t\}_{t \in \ZZ}\)
are satisfied by Gaussian functional FARIMA models,
including models with a common memory parameter
\cite{LiRobinsonShang2020} and models with coordinate-dependent
memory parameters \cite{RuizMedina2022,DurandRoueff2024},
under suitable summability and uniform spectral bounds. The simultaneous diagonalization condition also holds for fields in statistical physics whose spatial domain is a connected and compact two-point homogeneous space, provided that the covariance and spectral density operators are invariant under the isometry group of the manifold; see, for example, \cite{OvalleRM24}.

The proof relies on a Wiener--It\^o chaos decomposition for Hilbert-valued random variables (see Sections~2 and~4 of \cite{DukerZoub26}) and is split into two parts. First, we prove a reduction principle (cf.\ Proposition~\ref{prop:operator-reduction-principle}). In each basis direction \(v_i\) of \(\mbbH_2\), the fluctuations are governed by the leading terms, which share an effective memory exponent determined by their Hermite multi-indices and the memory parameters of the underlying Gaussian sequence.
For the proof of the reduction principle, we use finite Hermite projections and
uniform bounds on the projection remainders. Second, we prove that the kernels of the leading terms converge to the corresponding kernels of the limiting process \eqref{eq:ZD-wiener-chaos} (cf.\ Proposition~\ref{prop:operator-kernel-convergence}). Each individual term in this representation belongs to a fixed Wiener chaos. We establish the kernel convergence for each fixed Hermite multi-index and extend this convergence to the
full expansion by dominated convergence.
Moreover, the limiting process need not belong to a single
Wiener chaos, since different chaos orders may survive within a single
direction of $\mathbb{H}_2$, as well as across different directions.

The paper is organized as follows. Section~\ref{ftsm} introduces the spectral setting and some preliminaries for Hilbert-valued random variables. In Section~\ref{nes}, we state and discuss our assumptions and present the main operator-normalized theorem. Section~\ref{sec:hermite-processes} constructs the limiting process appearing in our theorem. Section~\ref{sec:proofs} contains the kernel convergence result and the reduction principle, together with their proofs and the proof of the main theorem. Section \ref{sec:auxiliary-lemmas} includes some auxiliary lemmas. Finally, Section~\ref{sec:examples} gives two applications of our result.

\section{Preliminaries}
\label{ftsm}
\label{preliminaries}

In this section we present some preliminary facts about random variables taking values in Hilbert spaces, as well as their Hermite and Wiener chaos expansions.

\subsection{Random variables on Hilbert spaces}
Let  $\mathbb{H}$ be a real separable Hilbert space with the inner product $\left\langle\cdot,\cdot\right\rangle_{\mathbb{H}},$ and denote by $\widetilde{\mathbb{H}}=\mathbb{H}+i\mathbb{H},$ its complex version,
 whose elements are of the form $ \psi=\psi_{1}+ i\psi_{2}$ for $\psi_{i}\in \mathbb{H}$, $i =1 ,2$, equipped with the inner product given by
\begin{equation}
\left\langle\psi_{1}+ i\psi_{2},\widetilde{\psi}_{1}+i\widetilde{\psi}_{2}\right\rangle_{\widetilde{\mathbb{H}}}=\left\langle \psi_{1},\widetilde{\psi}_{1}\right\rangle_{\mathbb{H}}
+\left\langle \psi_{2},\widetilde{\psi}_{2}\right\rangle_{\mathbb{H}}+i\left(\left\langle \psi_{2},\widetilde{\psi}_{1}\right\rangle_{\mathbb{H}}-\left\langle \psi_{1},\widetilde{\psi}_{2}\right\rangle_{\mathbb{H}}\right).
\label{ip}
\end{equation}

We denote by $(\Omega,\mathcal{F},P)$ the underlying probability space, and by $L^{2}(\Omega,\mathcal{F},P: \mbbH)$ the space of square-integrable random variables taking values in $\mbbH$, where $\mbbH$ may denote a real or complex Hilbert space. For compact operators on \(\mathbb{H}\) (where again $\mbbH$ may be real or complexified), fix an orthonormal basis $\{\psi_i\}_{in \in \NN}$ of $\mbbH$; let \(\Sone(\mathbb{H})\),  \(\Stwo(\mathbb{H})\), and $\mathcal{L}({\mathbb{H}})$ denote the trace-class, Hilbert--Schmidt, and bounded linear operators on $\mathbb{H}$, respectively. For an operator $\mathcal{A}$ with adjoint  $\mathcal{A}^{\star}$ the corresponding norms are given by 
\begin{equation}
    \begin{gathered}
\left\|\mathcal{A}\right\|_{\mathcal{S}_{1}(\mbbH)}
=
\sum_{k\geq 1}\left\langle \left[\mathcal{A}^{\star}
\mathcal{A}\right]^{1/2}(\psi_{k}),\psi_{k}
\right\rangle_{\mathbb{H}},
\\
\left\|\mathcal{A}\right\|_{\mathcal{S}_{2}(\mbbH)} 
= \sqrt{
\left\|\mathcal{A}^{\star}\mathcal{A}\right\|_{\mathcal{S}_{1}(\mathbb{H})}}, 
\qquad
\|\mathcal{A}\|_{\mathcal{L}(\mbbH)} =
 \sup_{\psi\in \mbbH;\ \|\psi\|_{\mathbb H}=1}
 \|\mathcal{A}(\psi)\|_{\mathbb{H}},
\label{HSN}        
    \end{gathered}
\end{equation}

Recall that \(\{X_t\}_{t\in\ZZ}\) denotes a centered, stationary, Gaussian, and \(\mbbH_1\)-valued random variables, with covariance operator $Q$, and denote the law of \(X_0\) by \(\gamma_Q\). Its autocovariance operators are
\begin{equation}
    \label{eq:covariance-operators}
    \Gamma(v)
    =
    \E[X_v\otimes X_0],
    \qquad
    \Gamma(v)h
    =
    \E\big[\langle X_0,h\rangle_{\mathbb{H}_1}X_v\big],
    \qquad v\in\ZZ,
\end{equation}
such that \(Q=\Gamma(0)\). From the spectral theorem, there is a basis $\{\phi_r\}_{r \geq 1}$ of $\mathbb{H}_1$ such that $ Q\phi_r=\lambda_r \phi_r$ for all $r \ge 1$, where $\sum_{r\ge1}\lambda_r=\Tr(Q)$. Without loss of generality, we can assume that $Q$ is nondegenerate (meaning $\langle Q h, h \rangle_{\mbbH_1} > 0$ for all $h \neq 0$) by excluding the coordinates $\phi_r$ with
\(\lambda_r=0\) from $\mathbb{H}_1$.

Following \cite{PanaretosTavakoli2013} and \cite{RuizMedina2022}, we define by $\mathscr F_\omega$ the spectral density operator of the Hilbert space-valued process $\{X_t\}_{ t \in \ZZ}$. Assume that $\mathscr F_\omega$ is a strongly measurable positive
self-adjoint trace-class operator for almost every $\omega\in[-\pi,\pi]$. Under Assumption \ref{ass:trace-admissibility} below, the autocovariance operators are given by
\begin{equation}
\Gamma(t)=\int_{-\pi}^{\pi}e^{it\omega}
                  \mathscr F_\omega\,d\omega, \qquad t \in \ZZ, \label{eq:covariance-from-density}
\end{equation}
where these integrals are understood in the Bochner sense and imply that $\Gamma(0)\in\Sone(\mathbb{H}_1)$, $\E\|X_0\|_{\mathbb{H}_1}^2<\infty$, and that 
\begin{equation}
    \label{eq:trace-admissibility-consequences}
   \|\Gamma(t)\|_{\Sone(\mathbb{H}_1)}
    \le
    \int_{-\pi}^{\pi}\Tr(\mathscr F_\omega)\,\dd\omega.
\end{equation}
Moreover, the Banach-valued Riemann--Lebesgue lemma gives $ \|\Gamma(t)\|_{\Sone(\mathbb{H}_1)}\to0$ as $|t| \to \infty$.
Under Assumption \ref{ass:trace-admissibility}, there is a common basis $\{\phi_j\}_{j \geq 1}$ diagonalizing the spectral operator family $\{\mathscr F_\omega\}_{\omega \in [-\pi,\pi]}$ for almost all $\omega \in [-\pi,\pi]$, admitting the family of eigenvalues $\{f_k(\omega)\}_{k \in \NN, \omega \in [-\pi,\pi]}$. As a result of Assumptions \ref{ass:trace-admissibility}\ref{ass:trace-mathF}--\ref{ass:diagonal-density}, the same basis diagonalizes $\{\Gamma(t)\}_{t \in \mathbb{Z}}$, the eigenvalues $\{\lambda_r\}_r$ of $Q$ are related to $\{f_k(\omega)\}_{k \in \NN, \omega \in [-\pi,\pi]}$ by  $ \lambda_r
    =\int_{-\pi}^{\pi} f_r(\omega)\,\dd\omega$, and
    we obtain
\begin{equation}
    \label{eq:diagonal-trace-sum}
    \sum_{k\ge1}\int_{-\pi}^{\pi}f_k(\omega)\,\dd\omega<\infty.
\end{equation}
Finally, we define the normalized spectral eigenvalues by
\begin{equation}
    \label{eq:sigma-and-coordinate-density}
    \widetilde f_r(\omega)=\frac{f_r(\omega)}{\lambda_r}, \quad \omega \in [-\pi,\pi], \qquad \text{so that } \int_{-\pi}^{\pi}  \widetilde f_r(\omega) \,\dd\omega = 1.
\end{equation}

\subsection{The infinite-dimensional Hermite expansion}
\label{idns}

In this section, we collect some facts about the white noise and Hermite expansions for random variables taking values in $\mbbH_1$. We mostly use the exposition of \cite{DukerZoub26}.

Define the white-noise mapping by, for $x,h \in \mbbH_1$,
\begin{equation}
\label{eq:white-noise-mapping}
W_h(x)
=
\sum_{r\ge1}\lambda_r^{-1/2}
\langle x,\phi_r\rangle_{\mathbb H_1}
\langle h,\phi_r\rangle_{\mathbb H_1}, \quad \text{ implying }  \quad
W_{\phi_r}(X_t)
=
\frac{\langle X_t,\phi_r\rangle_{\mathbb H_1}}{\sqrt{\lambda_r}},
\end{equation}
where the series converges in \(L^2(\mathbb H_1,\gamma_Q)\).

Let \(G:\mathbb H_1\to\mathbb H_2\) be centered and in $L^2(\mathbb H_1,\gamma_Q;\mathbb H_2)$, i.e., $\E\|G[X_0]\|_{\mathbb H_2}^2<\infty$. Denote with
\(\{v_i\}_{i\ge1}\) an orthonormal basis of \(\mathbb H_2\).

The probabilists' Hermite polynomials are defined by $H_0(z)=1$ and, for $p \ge 1$, $H_p(z)=(-1)^p e^{z^2/2}\frac{d^p}{dz^p}e^{-z^2/2}$. For standard jointly Gaussian \(U,V\),
\(\E[H_p(U)H_m(V)]=p!(\E[UV])^p\delta_{pm}\).
Set
\begin{equation}
    \label{eq:hermite-multi-indices}
    \begin{gathered}
    \mathcal L
    \doteq \Big\{
    \ell=(\ell_m)_{m\ge1}\in\NN_0^\infty:
                  |\ell|\doteq \sum_{m\ge1}\ell_m<\infty\Big\},\quad
    \mathcal L_p \doteq \{\ell\in\mathcal L:|\ell|=p\},
    \quad \ell!\doteq \prod_{m\ge1}\ell_m!,\\
    M_0\doteq 0,
    \qquad M_\ell\doteq \max\{m:\ell_m\ge1\}\quad(\ell\ne0).
    \end{gathered}
\end{equation}
The generalized Hermite polynomials and coefficients are
\begin{equation}
    \label{eq:generalized-hermite-polynomials}
    H_\ell[x]
    \doteq \prod_{m=1}^{M_\ell}H_{\ell_m}(W_{\phi_m}(x)),
    \qquad \ell\in\mathcal L,
\end{equation}
with the empty product equal to one. Then, from Lemma 2.1 in \cite{DukerZoub26}, we have the Hermite decomposition
\begin{equation}
    \label{eq:full-hermite-expansion}
    G[x]
    =\sum_{i\ge1}\sum_{\ell\in\mathcal L}
       c_{i,\ell}H_\ell[x]\otimes v_i, \quad \text{where}  \quad    c_{i,\ell}
    \doteq \frac{1}{\ell!}
       \langle G,H_\ell\otimes v_i\rangle_
       {L^2(\mathbb H_1,\gamma_Q)\otimes\mathbb H_2},
\end{equation}
where we have used the identification
\(L^2(\mathbb H_1,\gamma_Q;\mathbb H_2)
 \cong L^2(\mathbb H_1,\gamma_Q)\otimes\mathbb H_2\).
The notation \(H_\ell[x]\otimes v_i\) is identified with the usual multiplication
\(H_\ell[x]v_i\).

\begin{lemma}
\label{lem:hermite-basis-detail}
The family
\(\{(\ell!)^{-1/2}H_\ell\otimes v_i:
\ell\in\mathcal L,\ i\ge1\}\) is an orthonormal basis of
\(L^2(\mathbb H_1,\gamma_Q;\mathbb H_2)\). In particular,
\eqref{eq:full-hermite-expansion} converges in that Hilbert space and $\E\|G[X_0]\|_{\mathbb H_2}^2
   =\sum_{i\ge1}\sum_{\ell\in\mathcal L}\ell!|c_{i,\ell}|^2$.
\end{lemma}

For nonzero centered \(G\), its Hermite rank and chaotic expansion are defined by
\begin{equation}
    \label{eq:hermite-rank}
    q=\rank(G)
    \doteq \min\{p\ge1:c_{i,\ell}\ne0
                     \text{ for some }i\ge1,\ \ell\in\mathcal L_p\},
\end{equation}
\begin{equation}
    \label{eq:chaotic-components}
     G=\sum_{p=q}^{\infty}G_p,
   \qquad G_p[x]
    \doteq \sum_{i\ge1}\sum_{\ell\in\mathcal L_p}
       c_{i,\ell}H_\ell[x]\otimes v_i, \; p \ge q,
\end{equation}
and taking $G_0=\cdots=G_{q-1}=0$, where \(G_p\) denotes the chaos projection. Moreover, by the orthogonality of the Hermite polynomials we get the identities
\begin{align}
 \E\|G[X_0]\|_{\mathbb{H}_2}^2
 &=\sum_{p=q}^{\infty}\sum_{i=1}^\infty \sum_{\boldsymbol\ell\in\mathcal L_p}
        \boldsymbol\ell!\,|c_{i,\boldsymbol\ell}|^2<\infty,
 \quad
 \E\|G_q[X_0]\|_{\mathbb{H}_2}^2
 =\sum_{i=1}^\infty \sum_{\boldsymbol\ell\in\mathcal L_q}
        \boldsymbol\ell!\,|c_{i,\boldsymbol\ell}|^2<\infty.
 \label{eq:ordinary-coefficient-summability}
\end{align}

\subsection{Wiener chaos expansions on Hilbert spaces}
\label{subse-Wiender}

Let \(\widehat W\) be a Hermitian Gaussian spectral measure on
\(\mathbb R\times\mathbb N\), with control measure \(dx\otimes\#\),
so that its coordinate measures \(\widehat W_j\) are independent and
\(\widehat W_j(-A)=\overline{\widehat W_j(A)}\).
Let \(\mathfrak H\) be the real Hilbert space of complex-valued
\(L^2(dx\otimes\#)\) functions satisfying
\(f(-x,j)=\overline{f(x,j)}\), with the real \(L^2\) inner product.

For \(p\ge1\), let \(I_p\) denote the multiple Wiener--It\^o integral
of order \(p\) with respect to \(\widehat W\).
Equip the symmetric tensor powers with the norm inherited from
the full tensor power, keeping the factorial explicitly in the
isometry. For \(f\in\mathfrak H^{\odot p}\) and
\(g\in\mathfrak H^{\odot q}\), we then have
\begin{equation}
\label{eq:wiener-isometry}
 \E[I_p(f)I_q(g)]
 =
 \begin{cases}
 p!\langle f,g\rangle_{\mathfrak H^{\otimes p}},&p=q,\\
 0,&p\ne q.
 \end{cases}
\end{equation}
We denote its Hilbert-valued counterpart by
\(\mathcal I_p=I_p\otimes I_{\mathbb H_2}\).

Since the standardized coordinate process
\(\{W_{\phi_j}(X_t)\}_{t\in\mathbb Z}\) has spectral density
\(\widetilde f_j\), we may work with a joint spectral realization
such that
\begin{equation}
\label{eq:spectral-representation-coordinate}
 W_{\phi_j}(X_t)
 =
 \int_{-\pi}^{\pi}
 e^{it\omega}\widetilde f_j(\omega)^{1/2}
 \widehat W_j(d\omega) = I_1(h_{j,t}), \quad h_{j,t} (\omega,r) \doteq e^{it\omega}\widetilde f_j(\omega)^{1/2} \mathds{1}_{[-\pi,\pi]}(\omega) \mathds{1}_{\{r=j\}}.
\end{equation}

Let \(\operatorname{Sym}_p\) denote the standard symmetrization
operator from \(\mathfrak H^{\otimes p}\) onto
\(\mathfrak H^{\odot p}\); see
\cite[Appendix~B.3]{NourdinPeccati2012}.
For \(\boldsymbol\ell\in\mathcal L_p\), let
\(k_1\le\cdots\le k_p\) list the coordinate indices, with each
\(j\) appearing exactly \(\ell_j\) times.
For a square-integrable kernel \(f\) satisfying
\(f(-x)=\overline{f(x)}\), define
\begin{equation}
\label{eq:labelled-lift}
 f_{\boldsymbol\ell}((x_1,j_1),\ldots,(x_p,j_p))
 =
 f(x_1,\ldots,x_p)
 \prod_{a=1}^p\mathds{1}_{\{j_a=k_a\}}, \quad 
 J_{\boldsymbol\ell}f
 =
 \operatorname{Sym}_p f_{\boldsymbol\ell}.
\end{equation}
If \(f\) is invariant under permutations of positions carrying
the same label, then
\begin{equation}
\label{eq:labelled-isometry}
 p!\|J_{\boldsymbol\ell}f\|_{\mathfrak H^{\otimes p}}^2
 =
 \boldsymbol\ell!\|f\|_{L^2(\mathbb R^p)}^2.
\end{equation}

From \eqref{eq:spectral-representation-coordinate}, the labeled kernels satisfy
\(\|h_{j,t}\|_{\mathfrak H}=1\) and
\(\langle h_{j,t},h_{k,t}\rangle_{\mathfrak H}=0\) for \(j\ne k\).
Thus \(H_{\ell_j}(W_{\phi_j}(X_t))
=I_{\ell_j}(h_{j,t}^{\otimes\ell_j})\), and all cross-coordinate
contractions vanish in the product formula.
Expanding the resulting tensor product and applying
\(J_{\boldsymbol\ell}\) therefore yields
\begin{equation}
\label{eq:hermite-spectral-integral}
 H_{\boldsymbol\ell}[X_t]
 =
 I_p\left(
 J_{\boldsymbol\ell}
 \left[
 e^{it\sum_{a=1}^p\omega_a}
 \prod_{a=1}^p
 \left(
 \widetilde f_{k_a}(\omega_a)^{1/2}
 \mathds{1}_{\{|\omega_a|\le\pi\}}
 \right)
 \right]
 \right),
\end{equation}
where the spectral densities are extended by zero outside
\([-\pi,\pi]\).

\section{Assumptions and main results}
   \label{nes}

In this section, we present our assumptions, main results, and a discussion on the assumptions.

\subsection{Assumptions}

We start with the central assumption of our work. We call an input coordinate \(j\) \textit{active} if it appears in a nonzero term
of the Hermite expansion of \(G\), i.e., there is a multi-index $\ell$ and an $i$ such that $\ell_j > 0$ and $c_{i,\ell} \neq 0$.

We state the spectral assumptions on the underlying process $\{X_t\}_{ t\in \ZZ}$. 

\begin{assumption}
\label{ass:trace-admissibility}
Assume that the following hold for the spectral density operators $\{\mathscr F_\omega\}_{\omega}$:
\begin{enumerate}[label=(\roman*), ref=(\roman*), align=left]
\item \label{ass:trace-mathF}
For almost every \(\omega\in[-\pi,\pi]\), $\mathscr F_\omega\in\Sone(\widetilde{\mathbb{H}}_1)$ and $\int_{-\pi}^{\pi}\Tr(\mathscr F_\omega)\,\dd\omega  <\infty$.
\item \label{ass:diagonal-density}
There is an orthonormal basis \(\{\phi_k\}_{k\ge1}\) of \(\mathbb{H}_1\) such that, for all $k\ge1$,
\begin{equation}
    \label{eq:diagonal-density}
    \mathscr F_\omega \phi_k
    =
    f_k(\omega)\phi_k, 
    \quad
    f_k(\omega)\ge0,\quad f_k(-\omega)=f_k(\omega),
\end{equation}
for almost every \(\omega\in[-\pi,\pi]\). 

\item  \label{ass:coordinatewise-lrd}
There exist constants \(\bar\alpha\in(0,1)\) and \(0\le M<\infty\) such that, for every \(k\in\mathbb N\), there exist constants \(\alpha_k\in(0,\bar\alpha]\), \(\kappa_k>0\), and a measurable function \(L_k:(0,\pi]\to[0,\infty)\) satisfying
\begin{equation}
    \label{eq:global-spectral-envelope}
    f_k(\omega)
    =
    \lambda_k \kappa_k |\omega|^{-\alpha_k}L_k(|\omega|),
    \qquad
    \text{for a.e. }0<|\omega|\le\pi,
\end{equation}
where $\lim_{r\downarrow0}L_k(r)=1$ and $\sup_{k} L_k(r) \le M$ a.e. $0 < r \le |\pi|$.
\item \label{ass:uniform-spectral-scale}
Assume that the standardized spectral scales are uniformly bounded:
\begin{equation}
\label{eq:uniform-spectral-scale}
\kappa_*\doteq \sup\{ \kappa_k: k \; \text{is active}\}<\infty.
\end{equation}
\end{enumerate}
\end{assumption}

We state an important consequence of Assumption \ref{ass:trace-admissibility}.
Recall the memory parameters \(\{\alpha_j\}_{j \in \NN}\) from 
Assumption~\ref{ass:trace-admissibility}\ref{ass:coordinatewise-lrd}. For
\(\boldsymbol\ell\in\mathcal L_p\), define
\begin{equation}
\label{eq:unequal-effective-exponent}
 \beta(\boldsymbol\ell)
 =\sum_{j \in \supp(\boldsymbol\ell)} \ell_j(1-\alpha_j),
 \qquad
 H(\boldsymbol\ell)
 =1-\tfrac12\beta(\boldsymbol\ell).
\end{equation}
The quantity \(\beta(\boldsymbol\ell)\) combines the coordinate
memory parameters with their Hermite multiplicities, and we
call \(H(\boldsymbol\ell)\) the \textit{effective exponent}
of the corresponding Hermite term. The effective exponent corresponds to the scaling of our sums. We call the index $H(\ell)$ \textit{non-central} when $H(\ell) > 1/2$ or, equivalently, $\beta(\ell) < 1$. Denoting \(\varepsilon=1-\sup_j\alpha_j>0\), this implies that $\beta(\boldsymbol\ell)
 \ge\varepsilon\sum_j\ell_j
 =p\varepsilon$, hence
\(p<\varepsilon^{-1}\). In other words, for each fixed $\varepsilon$, we can only obtain Hermite polynomials of order up to $p < \varepsilon^{-1}$.

Define the dominating index and sets respectively by
\begin{equation}
\label{eq:dominant-exponent-and-set}
 \begin{gathered}
 d_i=\sup\{H(\boldsymbol\ell):c_{i,\boldsymbol\ell}\ne0,
 \ H(\boldsymbol\ell)>1/2\},\\
 \mathcal D_i
 =\{\boldsymbol\ell\in\mathcal L:c_{i,\boldsymbol\ell}\ne0,
\ H(\boldsymbol\ell)=d_i\}, \quad \mathcal N_i=\{\boldsymbol\ell\in\mathcal L:c_{i,\boldsymbol\ell}\ne0\}
                         \setminus\mathcal D_i,
 \end{gathered}
\end{equation}
giving rise to the dominating part of the Hermite expansion 
\begin{equation}
\label{eq:dominant-block}
 G_{\mathcal D}[X_t]
 =\sum_{i=1}^\infty \sum_{\boldsymbol\ell\in\mathcal D_i}
c_{i,\boldsymbol\ell}H_{\boldsymbol\ell}[X_t]v_i.
\end{equation}

\begin{assumption}
\label{ass:bigAssump2}
Fix a measurable, centered map $G:\mbbH_1 \to \mbbH_2$ with coefficients as in
\eqref{eq:full-hermite-expansion} and assume the following:
\begin{enumerate}[label=(\roman*), ref=(\roman*), align=left]
\item \label{ass:operator-standing} 
For every nonzero output coordinate, assume that its set of
active non-central indices is nonempty and that the supremum
defining \(d_i\) in \eqref{eq:dominant-exponent-and-set} is attained. For an identically zero output
coordinate, set \(d_i=3/4\) and \(\mathcal D_i=\varnothing\).

\item \label{ass:summability}
The following is satisfied
\begin{equation}
\label{eq:reduction-summability}
 \sum_{i=1}^{\infty}
 \frac{\E|\langle G[X_0],v_i\rangle_{\mathbb H_2}|^2}{2d_i-1}
 <\infty.
\end{equation}
\end{enumerate}
\end{assumption}

An important corollary of our main setting is the one where there is a common active component. In this simplified setting, Assumption \ref{ass:bigAssump2} will be replaced by the following, simplified hypothesis.
\begin{assumption}
\label{ass:common-active-exponent}
Suppose the following three assumptions hold:
\begin{enumerate}[label=(\roman*), ref=(\roman*), align=left]
\item \label{ass:GinL2}
The map $G \in L^2(\mbbH_1,\gamma_Q; \mbbH_2)$ has Hermite rank $q$.
\item \label{ass:common-active-dom-chaos} There is \(\alpha\in(0,1)\) such that
$c_{i,\boldsymbol\ell}\ne0$ for $\boldsymbol\ell\in\mathcal L_q$ implies 
$\alpha_j=\alpha$ for all $j\in\supp(\boldsymbol\ell)$, 
and $\alpha>1-\frac1q$. 
\item \label{ass:common-active-high-chaos} All higher-order chaoses scale more slowly than the leading chaos, i.e., if $c_{i,\boldsymbol\ell}\ne0$ and $|\boldsymbol\ell|>q$, then $H(\boldsymbol\ell)<1-\frac q2(1-\alpha)$.
\end{enumerate} 
\end{assumption}
For the common $\alpha$ in Assumption \ref{ass:common-active-exponent}, we define
\begin{equation}
\label{eq:common-H}
 H=1-\frac q2(1-\alpha)\in(1/2,1).
\end{equation}
Note that, if all active input coordinates, including those in higher chaoses, have
the common exponent $\alpha$, then Assumption \ref{ass:common-active-exponent}\ref{ass:common-active-high-chaos} follows immediately from Assumption \ref{ass:common-active-exponent}\ref{ass:common-active-dom-chaos}, since $1-p(1-\alpha)/2<H$ for $p>q$.

\subsection{Main theorems}
Recall the objects in \eqref{eq:dominant-exponent-and-set}.
Define the bounded diagonal operator and its normalization by
\begin{equation}
\label{eq:operator-normalization-D}
 Dv_i=d_iv_i,\qquad n^{-D}v_i=n^{-d_i}v_i
\end{equation}
and the partial-sum process and its continuous linear interpolation:
\begin{align}
S_n(s)=T_n\sum_{t=1}^{\lfloor ns\rfloor}G[X_t],
\qquad
\widetilde S_n(s)
 =S_n(s)+(ns-\lfloor ns\rfloor) T_n G[X_{\lfloor ns\rfloor+1}],
 \label{eq:continuous-interpolation}
\end{align}
where \(\{T_n\}_{n \geq 1}\) is a suitable sequence of scaling operators.
The random variable $\widetilde{S}_n$ is an element of
\(C([0,1];\mathbb{H}_2)\). The following is our main result, identifying the scaling limits for \eqref{eq:continuous-interpolation}.

\begin{theorem}
\label{thm:operator-normalized-full}
Let $\{X_t\}_{ t \in \ZZ}$ be a centered stationary Gaussian sequence
with values in $\mathbb H_1$, and let
$G\in L^2(\mathbb H_1,\gamma_{Q};\mathbb H_2)$ be centered
and nonzero, where
$\gamma_{Q}$ denotes the law of $X_0$. Suppose that
Assumptions~\ref{ass:trace-admissibility} and \ref{ass:bigAssump2}
hold. Let $D$ be the bounded diagonal operator in
\eqref{eq:operator-normalization-D}, and denote by $Z_D$ the
$\mathbb H_2$-valued process given in \eqref{eq:operator-limit-process} below.
Then, as $n\to\infty$,
\begin{equation}
\label{eq:operator-normalized-convergence}
 n^{-D}\sum_{t=1}^{\lfloor ns\rfloor}G[X_t]
 \xrightarrow{\operatorname{f.d.d.}} Z_D(s)
\end{equation}
in \(\mathbb{H}_2\), where $\xrightarrow{\operatorname{f.d.d.}}$ denotes convergence in finite-dimensional distribution. If, in addition, $d_*\doteq \inf_{i:\,\langle G,v_i\rangle_{\mathbb H_2}\ne0}d_i>\frac12$, then the process \(\widetilde S_n \) in
\eqref{eq:continuous-interpolation} with $T_n = n^{-D}$ converges weakly in
\(C([0,1];\mathbb{H}_2)\) to a continuous version of \(Z_D\). 
\end{theorem}

The limit $Z_D$ in \eqref{eq:operator-normalized-convergence} may contain several Hermite orders. If all surviving terms have
order one, the limit is Gaussian.

We consider the case where all input coordinates appearing in the
leading chaos have a common memory exponent, which is a corollary of Theorem~\ref{thm:operator-normalized-full}. 

\Needspace{18\baselineskip}
\begin{corollary}
\label{thm:common-hermite-limit}
Let $\{X_t\}_{ t \in \ZZ}$ be a centered stationary Gaussian sequence
with values in $\mathbb H_1$, and let
$G\in L^2(\mathbb H_1,\gamma_{Q};\mathbb H_2)$ be centered
and nonzero, where
$\gamma_{Q}$ denotes the law of $X_0$. Suppose that Assumptions~\ref{ass:trace-admissibility} and
\ref{ass:common-active-exponent} hold. Let $H$ and $Z_H^{(q)}$ be
given in \eqref{eq:common-H} and \eqref{eq:common-limit-process},
respectively. 
Then, as $n\to\infty$,
\begin{equation}
\label{eq:common-hermite-convergence}
 n^{-H}\sum_{t=1}^{\lfloor ns\rfloor}G[X_t]
 \xrightarrow{\operatorname{f.d.d.}} Z_H^{(q)}(s)
\end{equation}
in \(\mathbb{H}_2\) and, moreover,
the process \(\widetilde S_n \) in
\eqref{eq:continuous-interpolation} with $T_n = n^{-H} I_{\mathbb{H}_2}$ converges
weakly in \(C([0,1];\mathbb{H}_2)\) to a continuous version of \(Z_H^{(q)}\).
\end{corollary}

\subsection{Discussion on the assumptions}

\begin{remark}
\label{rem1}
As pointed out in the Introduction, an interesting example where the double
spectral diagonalization~\eqref{eq:diagonal-density} holds is obtained when the spatial
domain is a connected and compact two-point homogeneous space $\mathbb M_d$
with topological dimension $d$. Invariance under the group of isometries
leads to diagonalization of the covariance and spectral density operator
families in the eigenfunctions of the Laplace--Beltrami operator on
$\mathbb M_d$; see \cite{OvalleRM24}.
\end{remark}

\begin{remark}
\label{rem:diagonal-operator-formulation}
Replacing Assumption \ref{ass:trace-admissibility}\ref{ass:uniform-spectral-scale} with the slightly stronger assumption \(\sup_k\kappa_k<\infty\), we obtain a particularly convenient equivalent formulation of Assumption \ref{ass:trace-admissibility}\ref{ass:coordinatewise-lrd}. We can recast the spectral density operators as
\begin{equation}
 \mathscr F_\omega
 =
 Q^{1/2}|\omega|^{-A}L_\omega\mathcal K Q^{1/2},
 \qquad \text{for a.e. }0<|\omega|\le\pi,
\end{equation}
where the associated operators are diagonal and share a common eigenbasis 
\(Q\phi_k=\lambda_k\phi_k\), define
\[
 A\phi_k=\alpha_k\phi_k,\qquad
 \mathcal K\phi_k=\kappa_k\phi_k,\qquad
 L_\omega\phi_k=L_k(|\omega|)\phi_k.
\]
Then, we need that $ \|A\|_{\mathcal L(\mathbb H_1)}<1$, $\|\mathcal K\|_{\mathcal L(\mathbb H_1)}
 =\sup_k\kappa_k<\infty$, $0\le L_\omega\le M$ for a.e. $0<|\omega|\le\pi$, and that $L_\omega h \to h$ for every $h \in \mbbH_1$ as $\omega \to 0$.
\end{remark}

The following remark gives some sufficient time-domain assumptions ensuring our spectral Assumption~\ref{ass:trace-admissibility}.
\begin{remark} \label{rem:time-domain-sufficient}
Assumption~\ref{ass:trace-admissibility} holds
if \(Q\phi_k=\lambda_k\phi_k\), \(\lambda_k>0\),
\(\sum_k\lambda_k<\infty\), and the autocovariance operators
are diagonalized by the same basis and satisfy
\[
 \Gamma(h)\phi_k
 =
 \lambda_k c_k h^{\alpha_k-1}L_k(h)\phi_k,
 \qquad h\ge1,
\]
where
\(0<\inf_k\alpha_k\le\alpha_k\le\sup_k\alpha_k<1\),
\(0<\inf_k c_k\le c_k\le\sup_k c_k<\infty\), and each
\(L_k\) is measurable and slowly varying at infinity, with
\(\lim_{x\to\infty}L_k(x)=1\). We additionally require
\[
 \sup_k\sum_{h=1}^{\infty}
 h^{\alpha_k-1}|L_k(h)-1|<\infty.
\]

Indeed, the power-law cosine-series expansion, uniformly over
the stated compact range of memory parameters, together with
the weighted summability condition yields
\[
 \frac{f_k(\omega)}{\lambda_k}
 =
 \kappa_k|\omega|^{-\alpha_k}+R_k(\omega),
 \qquad
 \sup_k\|R_k\|_\infty<\infty,
 \qquad
 \kappa_k
 =
 \frac{c_k\Gamma(\alpha_k)\cos(\pi\alpha_k/2)}{\pi}.
\]
The expression for \(\kappa_k\) is the standard scalar
asymptotic constant; see
\cite[Proposition~2.2.14]{PipirasTaqqu17}.
Since the \(\kappa_k\)'s are uniformly bounded above and away
from zero, and \(\sup_k\|R_k\|_\infty<\infty\),
Assumption~\ref{ass:trace-admissibility} follows.
\end{remark}

\begin{remark}
Suppose that $G$ has Hermite rank $q$. Then, for all $\ell \in \mathcal{L}_p$, $p \ge q$, we have that $H(\boldsymbol \ell) = 1-\frac12\sum_{j\ge1}\ell_j(1-\alpha_j)
\le
1-\frac p2(1-\bar\alpha)
\le
1-\frac q2(1-\bar\alpha) $, and the equality in the final bound holds if and only if $p=q$ and $\alpha_j=\bar\alpha$ for every $j\in\operatorname{supp}(\boldsymbol\ell)$.
\end{remark}

\section{Construction of the limit}
\label{sec:hermite-processes}

In this section, we construct the limiting random variable $Z_D$ in \eqref{eq:operator-normalized-convergence}.

\subsection{Spectral kernels}
\label{subsec:kernels}

For \(s\in[0,1]\), define the temporal kernel
\begin{equation}
\label{eq:temporal-kernel}
 \widehat K_s(y)
 =
 \int_0^s e^{ivy}\,dv
 =
 \begin{cases}
 \dfrac{e^{isy}-1}{iy},&y\ne0,\\[1ex]
 s,&y=0.
 \end{cases}
\end{equation}
Fix \(\boldsymbol\ell\in\mathcal L_p\), \(p\ge1\), with
\(H(\boldsymbol\ell)>1/2\), and let \(k_1\le\cdots\le k_p\)
list the coordinate indices, with each \(j\) appearing exactly
\(\ell_j\) times. For \(x=(x_1,\ldots,x_p)\), set
\begin{align}
 K_{s,\boldsymbol\ell}(x)
 &=
 \widehat K_s\!\left(\sum_{a=1}^p x_a\right)
 \prod_{a=1}^p|x_a|^{-\alpha_{k_a}/2},
 \label{eq:operator-limit-kernel}\\
 L_{n,s,\boldsymbol\ell}(x)
 &=
 n^{-H(\boldsymbol\ell)-p/2}
 \sum_{t=1}^{\lfloor ns\rfloor}
 e^{it\sum_{a=1}^p x_a/n}
 \prod_{a=1}^p
 \left(
 \widetilde f_{k_a}(x_a/n)^{1/2}
 \mathds{1}_{\{|x_a|\le\pi n\}}
 \right).
 \label{eq:operator-prelimit-kernel}
\end{align}
Define, moreover,
\begin{equation}
\label{eq:kernel-envelope}
 B_{\boldsymbol\ell}
 =
 \int_{\mathbb R^p}
 (1\wedge|x_1+\cdots+x_p|^{-1})^2
 \prod_{a=1}^p\kappa_{k_a}|x_a|^{-\alpha_{k_a}}\,dx
\end{equation}
and, for \(p\ge2\), the convolution
\begin{equation} \label{eq:psi-form}
  \Psi_{\boldsymbol\ell}(y)
 =
 \int_{\mathbb R^{p-1}}
 \left(
 \prod_{a=1}^{p-1}
 \kappa_{k_a}|x_a|^{-\alpha_{k_a}}
 \right)
 \kappa_{k_p}
 \left|y-\sum_{a=1}^{p-1}x_a\right|^{-\alpha_{k_p}}
 \,dx_1\cdots dx_{p-1},   
\end{equation}
while, for \(p=1\), set
\(\Psi_{\boldsymbol\ell}(y)
=\kappa_{k_1}|y|^{-\alpha_{k_1}}\).

\begin{lemma}
\label{lem:envelope-finiteness}
For every \(\boldsymbol\ell\in\mathcal L_p\), \(p\ge1\), satisfying
\(H(\boldsymbol\ell)>1/2\), we have
\(0<\Psi_{\boldsymbol\ell}(1)<\infty\) and
\begin{equation}
\label{eq:power-convolution}
 \Psi_{\boldsymbol\ell}(y)
 =
 \Psi_{\boldsymbol\ell}(1)
 |y|^{1-2H(\boldsymbol\ell)},
 \qquad y\ne0.
\end{equation}
Moreover,
\begin{equation}
\label{eq:envelope-convolution-integral}
 B_{\boldsymbol\ell}
 =
 \Psi_{\boldsymbol\ell}(1)
 \left(
 \frac1{1-H(\boldsymbol\ell)}
 +
 \frac1{H(\boldsymbol\ell)}
 \right)
 <\infty.
\end{equation}
\end{lemma}

\begin{proof}
For \(p=1\), we have
\(1-2H(\boldsymbol\ell)=-\alpha_{k_1}\) and
\(\Psi_{\boldsymbol\ell}(1)=\kappa_{k_1}\), so
\eqref{eq:power-convolution} follows directly from the
definition of \(\Psi_{\boldsymbol\ell}\).

For \(p\ge2\) and $1\le r\le p$, set
\[
 \eta_r
 =
 \sum_{a=1}^r\alpha_{k_a}-r+1
 =
 1-\sum_{a=1}^r(1-\alpha_{k_a}).
\]
Since the terms \(1-\alpha_{k_a}\) are all positive and their
total is \(\beta(\boldsymbol\ell)<1\), \( \eta_r \in (0,1)\) for all $1 \le r \le p$.

We use the form of $\Psi$ in \eqref{eq:psi-form} to prove \eqref{eq:power-convolution}. Note that it suffices to prove that $\Psi_\ell(y) = C \Psi_\ell(1) |y|^{1-2H(\ell)}$, since then we can obtain the constant by evaluating $\Psi_\ell$ at one.

We prove that, for every \(r=1,\ldots,p\), there exists
\(0<C_r<\infty\) such that, whenever $y\ne0$,
\[
 \int_{\mathbb R^{r-1}}
 \left(\prod_{a=1}^{r-1}
       \kappa_{k_a}|x_a|^{-\alpha_{k_a}}\right)
 \kappa_{k_r}
 \left|y-\sum_{a=1}^{r-1}x_a\right|^{-\alpha_{k_r}}
 \,dx_1\cdots dx_{r-1}
 =
 C_r|y|^{-\eta_r}.
\]
For \(r=1\), the left-hand side is understood as
\(\kappa_{k_1}|y|^{-\alpha_{k_1}}\).
Since \(\eta_1=\alpha_{k_1}\), the identity holds with
\(C_1=\kappa_{k_1}\).

Fix \(2\le r\le p\) and assume the identity holds for
\(r-1\). By Tonelli's theorem and the induction hypothesis,
\begin{multline*}
\int_{\mathbb R^{r-1}}
\left(\prod_{a=1}^{r-1}
\kappa_{k_a}|x_a|^{-\alpha_{k_a}}\right)
\kappa_{k_r}
\left|y-\sum_{a=1}^{r-1}x_a\right|^{-\alpha_{k_r}}
dx_1\cdots dx_{r-1} 
= 
C_{r-1}\kappa_{k_r}
\int_{\mathbb R}
|u|^{-\eta_{r-1}}|y-u|^{-\alpha_{k_r}} du 
\\
=
C_{r-1}\kappa_{k_r}|y|^{-\eta_r}
\int_{\mathbb R}
|v|^{-\eta_{r-1}}|1-v|^{-\alpha_{k_r}} dv = C_r |y|^{-\eta_r},
\end{multline*}
where for the second equality we used the substitution $u=yv$ and that 
\(\eta_r=\eta_{r-1}+\alpha_{k_r}-1\); for the third equality we have set $C_r
 =
 C_{r-1}\kappa_{k_r}
 \int_{\mathbb R}
 |v|^{-\eta_{r-1}}|1-v|^{-\alpha_{k_r}}\,dv
 \in(0,\infty)$, where we notice that the finiteness of the integral near zero and one is due to
\(0<\eta_{r-1},\alpha_{k_r}<1\), and at infinity because $\eta_{r-1}+\alpha_{k_r}=1+\eta_r>1$. This proves the induction step.

To obtain \(B_{\boldsymbol\ell}\), replace \(x_p\) by
\(y=x_1+\cdots+x_p\). Then, using Tonelli's theorem and the definition of
\(\Psi_{\boldsymbol\ell}\),
\begin{equation*}
     B_{\boldsymbol\ell}
 =
 \int_{\mathbb R}
 (1\wedge|y|^{-1})^2\Psi_{\boldsymbol\ell}(y)\,dy
 =
 2\Psi_{\boldsymbol\ell}(1)
 \left(
 \int_0^1 y^{1-2H(\boldsymbol\ell)}\,dy
 +
 \int_1^\infty y^{-1-2H(\boldsymbol\ell)}\,dy
 \right).
\end{equation*}
Both integrals are finite because
\(1/2<H(\boldsymbol\ell)<1\), proving \eqref{eq:envelope-convolution-integral}.
\end{proof}

\subsection{Definition of the limit}
\label{subsec:def-limit}

\subsubsection{The general case}

On the common spectral probability
space introduced in Section~\ref{subse-Wiender}, define, with $\kappa^\ell=\prod_j\kappa_j^{\ell_j}$,
\begin{align}
Z_D(s)=\sum_{i=1}^\infty  Z_i(s)v_i
\quad
\text{ with }
\quad
Z_i(s)=\sum_{\boldsymbol\ell\in\mathcal D_i}
c_{i,\boldsymbol\ell}\kappa^{\boldsymbol\ell/2}
I_{|\boldsymbol\ell|}(J_{\boldsymbol\ell}K_{s,\boldsymbol\ell}),
 \label{eq:operator-limit-process}
\end{align}
or, equivalently, grouping the terms of order \(p\) into a kernel
\(h_{p,s}\in\mathfrak H^{\odot p}\otimes\mathbb{H}_2\), we can write
\begin{equation} \label{eq:ZD-wiener-chaos}
Z_D(s)
=
\sum_{p\ge q}\mathcal I_p(h_{p,s}), 
\qquad 
h_{p,s} \doteq 
\sum_{i=1}^\infty
\sum_{\ell \in \mathcal D_i, |\boldsymbol \ell| = p}
c_{i,\boldsymbol\ell}\kappa^{\boldsymbol\ell/2}
J_{\boldsymbol\ell}K_{s,\boldsymbol\ell} \otimes v_i .
\end{equation}

By using \eqref{eq:labelled-isometry}, the fact that label supports for different multiplicity indices are disjoint, and the kernel estimates~\eqref{eq:uniform-kernel-norm} below, we can infer that
\begin{multline}
\label{eq:operator-limit-second-moment}
 \E\|Z_D(s)\|_{\mathbb{H}_2}^2
 =\sum_{i=1}^\infty \sum_{\boldsymbol\ell\in\mathcal D_i}
 \boldsymbol\ell!\,|c_{i,\boldsymbol\ell}|^2
 \|\kappa^{\boldsymbol\ell/2}K_{s,\boldsymbol\ell}\|_{L^2(\mathbb R^{|\boldsymbol\ell|})}^2 \\\le \sum_{i=1}^\infty \sum_{\boldsymbol\ell\in\mathcal D_i}
 \boldsymbol\ell!\,|c_{i,\boldsymbol\ell}|^2 C/(2d_i-1) \le C\sum_{i=1}^\infty \frac{\E|\langle G[X_0],v_i\rangle_{\mathbb H_2}|^2}{2d_i-1}<\infty,
\end{multline}
from Assumption \ref{ass:bigAssump2}\ref{ass:summability}.

\subsubsection{The common exponent case}
Let $q$ and $\alpha$ be as in Assumption \ref{ass:common-active-exponent}. Define
\[
K_s(x_1,\ldots,x_q)
\doteq
\widehat K_s\!\left(\sum_{a=1}^{q}x_a\right)
\prod_{a=1}^{q}|x_a|^{-\alpha/2},
\qquad (x_1,\ldots,x_q)\in\mathbb R^q,
\]
where $\widehat K_s$ is given by \eqref{eq:temporal-kernel}.
By Assumption \ref{ass:common-active-exponent}\ref{ass:common-active-dom-chaos} and \eqref{eq:operator-limit-kernel}, we have
$K_{s,\ell}=K_s$ for every $\ell\in\mathcal L_q$
such that $c_{i,\ell}\ne0$ for some $i$.
We then define
\begin{equation}
\label{eq:common-limit-process}
 Z_H^{(q)}(s)=\sum_{i=1}^\infty  I_q(\mathcal K_{i,s})v_i
 =\mathcal I_q\left(\sum_{i=1}^\infty \mathcal K_{i,s}\otimes v_i\right), \quad \mathcal K_{i,s}
 =\sum_{\boldsymbol\ell\in\mathcal L_q}
   c_{i,\boldsymbol\ell}\kappa^{\boldsymbol\ell/2}
   J_{\boldsymbol\ell}K_s
 \quad\text{in }\mathfrak H^{\odot q}.
\end{equation}
Note that when $s > 0$, $K_s\in L^2(\mathbb R^q)$ if and only if $1-1/q<\alpha<1$. The isometries \eqref{eq:wiener-isometry}--\eqref{eq:labelled-isometry} then imply
\begin{equation}
\label{eq:common-limit-second-moment}
 \E\|Z_H^{(q)}(s)\|_{\mathbb{H}_2}^2
 =\|K_s\|_{L^2(\mathbb R^q)}^2\sum_{i=1}^\infty \sum_{\boldsymbol\ell\in\mathcal L_q}
   \boldsymbol\ell!\,|c_{i,\boldsymbol\ell}|^2
   \kappa^{\boldsymbol\ell}<\infty.
\end{equation}

\section{Proofs of the limit results}
\label{sec:proofs}

We prove a reduction principle and derive convergence of the kernel functions that correspond to pre- and postlimits. Finally, we combine these two results to prove
Theorem~\ref{thm:operator-normalized-full}, and conclude with the
proof of Corollary~\ref{thm:common-hermite-limit}.

\subsection{Reduction principle}

For the following proofs, we first introduce the required notation.
Recall $\mathcal{L}$ from \eqref{eq:hermite-multi-indices}. Let $\mathcal Q=(\mathcal Q_i)_{i\ge1}$, where each
$\mathcal Q_i\subseteq\mathcal L\setminus\{0\}$ consists of indices supported on active coordinates. Fix $d_i\in(1/2,1)$ and recall the operator 
$D$ from \eqref{eq:operator-normalization-D}.
For $N\ge1$, define
\begin{equation}
\label{eq:def:set:mathcalF-N-D}
\mathcal F_N(\mathcal Q)
\doteq
\left\{(i,\boldsymbol\ell):
1\le i\le N,\quad \boldsymbol\ell\in\mathcal Q_i,\quad
|\boldsymbol\ell|\le N,\quad \ell_j=0\text{ for all }j>N\right\},
\end{equation}
and write
\begin{equation}
\label{eq:def:projections:G}
G_{\mathcal Q}
=\sum_{i=1}^{\infty}\sum_{\boldsymbol\ell\in\mathcal Q_i} c_{i,\boldsymbol\ell}H_{\boldsymbol\ell}v_i,
\qquad
G_{\mathcal Q}^{(N)}
=\sum_{(i,\boldsymbol\ell)\in\mathcal F_N(\mathcal Q)} c_{i,\boldsymbol\ell}H_{\boldsymbol\ell}v_i.
\end{equation}
The sets $\mathcal F_N(\mathcal Q)$ are finite and increase to all
indices in $\mathcal Q$. By Hermite orthogonality and Parseval's
identity, $G_{\mathcal Q}^{(N)}$ is an orthogonal projection of
$G_{\mathcal Q}$ onto a finite-dimensional subspace of
$L^2(\mathbb H_1,\gamma_Q;\mathbb H_2)$ and converges to
$G_{\mathcal Q}$ in this space. Each projection takes values in
$\operatorname{span}\{v_1,\ldots,v_N\}$  and depends only on the first
$N$ standardized Gaussian coordinates.

\begin{proposition}
\label{prop:operator-reduction-principle}
Recall $G_{\mathcal{D}}$ from 
\eqref{eq:dominant-block}.
Suppose Assumptions~\ref{ass:trace-admissibility}\ref{ass:diagonal-density}--\ref{ass:uniform-spectral-scale},
\ref{ass:bigAssump2}\ref{ass:operator-standing},\ref{ass:summability}. Then, as $n\to\infty$,
\begin{equation}
\label{eq:operator-reduction-result}
 \sup_{s\in[0,1]}\E\left\|
 n^{-D}\sum_{t=1}^{\lfloor ns\rfloor}
          (G[X_t]-G_{\mathcal D}[X_t])
 \right\|_{\mathbb H_2}^2\to0.
\end{equation}
\end{proposition}

\begin{proof}
We control each nondominant term uniformly over partial-sum
lengths. Recall $\mathcal N_i$ from \eqref{eq:dominant-exponent-and-set}, i.e., the set of indices that corresponds to nondominat terms.
Fix $i$ and $\boldsymbol\ell\in\mathcal N_i$.
By Assumption~\ref{ass:bigAssump2}\ref{ass:operator-standing} and
the definition of $\mathcal N_i$, we have
$H(\boldsymbol\ell)<d_i$ and $d_i>1/2$.
We may therefore choose $d'$ such that
$\max\{1/2,H(\boldsymbol\ell)\}<d'<d_i$.
Applying \eqref{eq:spectral-block-bound} of
Lemma~\ref{lem:spectral-variance} with $d=d'$ gives
\begin{align}
 \max_{0\le m\le n}n^{-2d_i}
 \Var\left(\sum_{t=1}^mH_{\boldsymbol\ell}[X_t]\right)
 &\le\frac{C\boldsymbol\ell!}{2d'-1}
 n^{-2(d_i-d')}
 \to0.
 \label{eq:reduction-termwise-bound-1}
\end{align}
For the inequality, we use $m^{2d'}\le n^{2d'}$ for
$1\le m\le n$; the case $m=0$ is immediate.
The convergence follows from $d'<d_i$.

Set $F=G-G_{\mathcal D}=G_{\mathcal N}$ and
$F^{(N)}=G_{\mathcal N}^{(N)}$, using the representation in \eqref{eq:def:projections:G} with $\mathcal{N} = (\mathcal{N}_i)_{i \geq 1}$.
For each fixed $N$, Lemma \ref{lem:hermite-basis-detail} gives
\begin{align}
 \sup_{s\in[0,1]}\E\left\|
 n^{-D}\sum_{t=1}^{\lfloor ns\rfloor}F^{(N)}[X_t]
 \right\|_{\mathbb H_2}^2
 &\le
 \sum_{(i,\boldsymbol\ell)\in\mathcal F_N(\mathcal N)}
 |c_{i,\boldsymbol\ell}|^2
 \max_{0\le m\le n}n^{-2d_i}
 \Var\left(\sum_{t=1}^mH_{\boldsymbol\ell}[X_t]\right)
 \to0.
 \label{eq:reduction-finite-projection}
\end{align}
The sum is finite and each term tends to zero by
\eqref{eq:reduction-termwise-bound-1}.
Recall that $H(\boldsymbol\ell)\le d_i$ for every
$\boldsymbol\ell\in\mathcal N_i$, so
Lemma~\ref{lem:finite-hermite-projection} applies with
$\mathcal Q=\mathcal N$. Combining its bound with
\eqref{eq:reduction-finite-projection} and the decomposition
$F=F^{(N)}+(F-F^{(N)})$, we obtain, for every $N\ge1$,
\begin{align*}
\limsup_{n\to\infty}\sup_{s\in[0,1]}\E\left\|
 n^{-D}\sum_{t=1}^{\lfloor ns\rfloor}F[X_t]
 \right\|_{\mathbb H_2}^2
 &\le
 2\limsup_{n\to\infty}\sup_{s\in[0,1]}\E\left\|
 n^{-D}\sum_{t=1}^{\lfloor ns\rfloor}F^{(N)}[X_t]
 \right\|_{\mathbb H_2}^2\\
 &\qquad\quad+
 2\limsup_{n\to\infty}\sup_{s\in[0,1]}\E\left\|
 n^{-D}\sum_{t=1}^{\lfloor ns\rfloor}
 (F-F^{(N)})[X_t]
 \right\|_{\mathbb H_2}^2\\
 &\le
 2C\sum_{\substack{i\ge1,\ \boldsymbol\ell\in\mathcal N_i\\
                   (i,\boldsymbol\ell)\notin\mathcal F_N(\mathcal N)}}
 \frac{\boldsymbol\ell!|c_{i,\boldsymbol\ell}|^2}{2d_i-1}\\
 &=
 2C\sum_{i=1}^{\infty}\sum_{\boldsymbol\ell\in\mathcal N_i}
 \frac{\boldsymbol\ell!|c_{i,\boldsymbol\ell}|^2}{2d_i-1}
 -2C\sum_{(i,\boldsymbol\ell)\in\mathcal F_N(\mathcal N)}
 \frac{\boldsymbol\ell!|c_{i,\boldsymbol\ell}|^2}{2d_i-1}
\to 0,
\end{align*}
as $N\to\infty$.
The first inequality follows from
$\|u+v\|_{\mathbb H_2}^2\le2\|u\|_{\mathbb H_2}^2+2\|v\|_{\mathbb H_2}^2$.
For the second inequality, the first term vanishes by
\eqref{eq:reduction-finite-projection}, while the second is bounded
by Lemma~\ref{lem:finite-hermite-projection}, uniformly in $n$ and $s$.
The convergence
follows, since, by Parseval's identity and
Assumption~\ref{ass:bigAssump2}\ref{ass:summability}, we have
\begin{align}
 \sum_{i=1}^{\infty}\sum_{\boldsymbol\ell\in\mathcal N_i}
 \frac{\boldsymbol\ell!|c_{i,\boldsymbol\ell}|^2}{2d_i-1}
 &\le\sum_{i=1}^{\infty}\frac{1}{2d_i-1}
 \sum_{\boldsymbol\ell\in\mathcal L}
 \boldsymbol\ell!|c_{i,\boldsymbol\ell}|^2
 =\sum_{i=1}^{\infty}
 \frac{\E|\langle G[X_0],v_i\rangle_{\mathbb H_2}|^2}{2d_i-1}
 <\infty.
 \label{eq:reduction-summable-majorant}
\end{align}
\end{proof}

\subsection{Kernel convergence}

Recall the prelimit kernels $L_{n,s,\boldsymbol\ell}$ from
\eqref{eq:operator-prelimit-kernel} and the limit kernels
$K_{s,\boldsymbol\ell}$ from~\eqref{eq:operator-limit-kernel}. The following proposition states a convergence result for these kernels.

\begin{proposition}
\label{prop:operator-kernel-convergence}
Suppose  Assumptions~\ref{ass:trace-admissibility}\ref{ass:diagonal-density}--\ref{ass:uniform-spectral-scale}.
For every \(\boldsymbol\ell\in\mathcal L_p\) supported on active
coordinates with \(H(\boldsymbol\ell)>1/2\), and every \(s\in[0,1]\), there is a constant $C$,
independent of the index and its order, such that
\begin{equation}
\label{eq:uniform-kernel-norm}
\sup_{n\ge1}\sup_{s\in[0,1]}\|L_{n,s,\boldsymbol\ell}\|_{L^2(\mathbb R^p)}^2
\le\frac{C}{2H(\boldsymbol\ell)-1},\qquad
\sup_{s\in[0,1]}\|\kappa^{\boldsymbol\ell/2}K_{s,\boldsymbol\ell}\|_{L^2(\mathbb R^p)}^2
 \le\frac{C}{2H(\boldsymbol\ell)-1},
\end{equation}
Moreover, for every fixed $s\in[0,1]$, as $n\to\infty$,
\[
L_{n,s,\boldsymbol\ell}
\to
\kappa^{\boldsymbol\ell/2}K_{s,\boldsymbol\ell}
\quad\text{in }L^2(\mathbb R^p).
\]
\end{proposition}

\begin{proof}
We first establish convergence in $L^2(\mathbb R^p)$ of the kernels
$L_{n,s,\boldsymbol\ell}$ to
$\kappa^{\boldsymbol\ell/2}K_{s,\boldsymbol\ell}$. This convergence is then used to infer \eqref{eq:uniform-kernel-norm}.

For the $L^2(\mathbb R^p)$-convergence, we first establish almost-everywhere convergence and then verify that dominated convergence theorem may be applied by showing integrability of the respective objects.
Fix
$\boldsymbol\ell\in\mathcal L_p$ with
$H(\boldsymbol\ell)>1/2$, write
$\mathbf k=\mathbf k(\boldsymbol\ell)$, and
\begin{equation}
\widehat K_{n,s}(y)\doteq \frac1n\sum_{t=1}^{\lfloor ns\rfloor}e^{ity/n}, \qquad
y=x_1+\cdots+x_p.
\label{eq:def-kernelhat}
\end{equation}
Fix also $s\in(0,1]$; the case $s=0$ is immediate because both
kernels vanish.
Using the definition in \eqref{eq:operator-prelimit-kernel}, we have, with further explanations given below,
\begin{align}
L_{n,s,\boldsymbol\ell}(x)
&=
n^{-H(\boldsymbol\ell)-p/2}
\sum_{t=1}^{\lfloor ns\rfloor}e^{it\sum_{a=1}^p x_a/n}
\prod_{a=1}^p\widetilde{f}_{k_a}(x_a/n)^{1/2}
\mathbf1_{\{|x_a|\le\pi n\}}
 \label{eq:kernel-factorization}
\\
&=\widehat K_{n,s}(y)\prod_{a=1}^p
 \left[n^{-\alpha_{k_a}/2}\widetilde f_{k_a}(x_a/n)^{1/2}
\mathbf1_{\{|x_a|\le\pi n\}}\right]
 \label{eq:kernel-factorization+1}
\\
&\to
 \widehat K_s(y)
 \prod_{a=1}^p
 \kappa_{k_a}^{1/2}|x_a|^{-\alpha_{k_a}/2}
=\kappa^{\boldsymbol\ell/2}K_{s,\boldsymbol\ell}(x)
 \quad\text{for a.e. }x\in\mathbb R^p,
  \label{eq:kernel-factorization+2}
\end{align}
where \eqref{eq:kernel-factorization+1} follows by \eqref{eq:def-kernelhat} and since
$-1-\frac12\sum_{a=1}^p\alpha_{k_a}
=-H(\boldsymbol\ell)-p/2$.
For \eqref{eq:kernel-factorization+2}, the convergence of each spectral factor follows from
Assumption~\ref{ass:trace-admissibility}\ref{ass:coordinatewise-lrd}, since $x_a/n\to0$ and,
for fixed $x_a\ne0$, the indicator equals one for all sufficiently
large $n$. The convergence $\widehat K_{n,s}(y) \to \widehat K_s(y)$, with $\widehat K_s$ defined in \eqref{eq:temporal-kernel},  follows from the end of p.\ 289, (5.3.10), (5.3.13), and p.\ 291 of
\cite{PipirasTaqqu17}.
The final equality uses
$\kappa^{\boldsymbol\ell/2} \doteq \prod_{a=1}^p
\kappa_{k_a}^{1/2}$ and the definition of $K_{s,\boldsymbol\ell}$ in \eqref{eq:operator-limit-kernel}.

For integrability, we introduce the set $A_n=\{x\in\mathbb R^p:|x_1+\cdots+x_p|\le\pi n\}$ and consider both $A_n$ and its complement separately.
On $A_n$, we get
\begin{align}
 |\widehat K_{n,s}(y)|
 &\le
 \min\left\{
 \frac{\lfloor ns\rfloor}{n},
 \frac1{n|\sin(y/(2n))|}
 \right\}
 \le
 \min\left\{
 1,\frac1{n|\sin(y/(2n))|}
 \right\}
 \\
 &\le \min\{1,\pi|y|^{-1}\}
 \le \pi(1\wedge|y|^{-1}),
 \qquad |y|\le\pi n,
 \label{eq:kernel-central-dirichlet}
\end{align}
where, in the first line, the two bounds follow respectively from the
triangle inequality and (A.2.2) in
\cite{PipirasTaqqu17}; we then use
$\lfloor ns\rfloor/n\le1$.
On $|y|\le\pi n$, the inequality
$|\sin(y/(2n))|\ge |y|/(\pi n)$ gives the second line. Then, using \eqref{eq:kernel-factorization+1},
\begin{align}
 |L_{n,s,\boldsymbol\ell}(x)\mathbf1_{A_n}(x)|^2
 &=
\left| \widehat K_{n,s}(y)\prod_{a=1}^p
 \left[n^{-\alpha_{k_a}/2}\widetilde f_{k_a}(x_a/n)^{1/2}
\mathbf1_{\{|x_a|\le\pi n\}}\right]  \right|^2 \mathbf1_{A_n}(x)
 \\
 &\le
M^p\kappa^{\boldsymbol\ell}
 |\widehat K_{n,s}(y)|^2
 \prod_{a=1}^p|x_a|^{-\alpha_{k_a}}
 \mathbf1_{A_n}(x)
\label{eq:kernel-global-domination}
 \\&
 \leq
\pi^2M^p\kappa^{\boldsymbol\ell}
 (1\wedge|y|^{-1})^2
 \prod_{a=1}^p|x_a|^{-\alpha_{k_a}}
 \doteq g_{\boldsymbol\ell}(x),
\label{eq:kernel-global-domination+1}
\end{align}
where \eqref{eq:kernel-global-domination} follows by Assumption \ref{ass:trace-admissibility}\ref{ass:coordinatewise-lrd} 
and $\prod_{a=1}^p\kappa_{k_a}=\kappa^{\boldsymbol\ell}$; \eqref{eq:kernel-global-domination+1} is due to 
\eqref{eq:kernel-central-dirichlet}. 
Combining this bound with
\eqref{eq:kernel-central-dirichlet} on $A_n$ gives the last line.
It remains to show that $g_{\boldsymbol\ell}(x)$ is integrable, that is,
\begin{align}
 \int_{\mathbb R^p}g_{\boldsymbol\ell}(x)\,dx
 &=\pi^2M^p
 \int_{\mathbb R^p}
 (1\wedge|x_1+\cdots+x_p|^{-1})^2
 \prod_{a=1}^p
 \kappa_{k_a}|x_a|^{-\alpha_{k_a}}\,dx
 =\pi^2M^pB_{\boldsymbol\ell}<\infty,
\end{align}
with $B_{\boldsymbol\ell}$ as in \eqref{eq:kernel-envelope}.
Finiteness was established in \eqref{eq:envelope-convolution-integral}.
By \eqref{eq:kernel-factorization+2}, for almost every $x\in\mathbb R^p$,
\begin{align}
L_{n,s,\boldsymbol\ell}(x)\mathbf1_{A_n}(x)
\to \kappa^{\boldsymbol\ell/2}K_{s,\boldsymbol\ell}(x),
\end{align}
using $\mathbf1_{A_n}(x)
\to1$ that holds since every fixed $x$ belongs to $A_n$ for all sufficiently large $n$.
Finally,
\begin{align}
 \left|
 L_{n,s,\boldsymbol\ell}(x)\mathbf1_{A_n}(x)
 -\kappa^{\boldsymbol\ell/2}K_{s,\boldsymbol\ell}(x)
 \right|^2
 &\le
 2|L_{n,s,\boldsymbol\ell}(x)\mathbf1_{A_n}(x)|^2
 +2|\kappa^{\boldsymbol\ell/2}K_{s,\boldsymbol\ell}(x)|^2
\le4g_{\boldsymbol\ell}(x),
\end{align}
using $|u-v|^2\le2|u|^2+2|v|^2$ and
since by \eqref{eq:kernel-factorization+2}, we have $|\kappa^{\boldsymbol\ell/2}K_{s,\boldsymbol\ell}(x)|^2
=\lim_{n\to\infty}
 |L_{n,s,\boldsymbol\ell}(x)\mathbf1_{A_n}(x)|^2
 \le g_{\boldsymbol\ell}(x)$.
The integrand in the last line converges to zero almost
everywhere and is bounded by the integrable function
$4g_{\boldsymbol\ell}$.

We can apply dominated convergence on
$A_n$ to prove \eqref{eq:kernel-central-convergence}
\begin{align}
 \left\|
 L_{n,s,\boldsymbol\ell}\mathbf1_{A_n}
 -\kappa^{\boldsymbol\ell/2}K_{s,\boldsymbol\ell}
 \right\|_{L^2(\mathbb R^p)}^2
 &=\int_{\mathbb R^p}
 \left|
 L_{n,s,\boldsymbol\ell}(x)\mathbf1_{A_n}(x)
 -\kappa^{\boldsymbol\ell/2}K_{s,\boldsymbol\ell}(x)
 \right|^2\,dx
 \to0.
 \label{eq:kernel-central-convergence}
\end{align}
It remains to control the prelimit kernel on $A_n^c$.
On $A_n^c$, use the integral of $|\widehat K_{n,s}|^2$. Parseval's identity in (A.1.4) of
\cite{PipirasTaqqu17} gives
\begin{equation}
 \int_{-\pi n}^{\pi n}|\widehat K_{n,s}(y)|^2\,dy
 =\frac1n\int_{-\pi}^{\pi}
   \left|\sum_{t=1}^{\lfloor ns\rfloor}e^{itx}\right|^2\,dx
 =\frac{2\pi\lfloor ns\rfloor}{n}
 \le2\pi.
 \label{eq:kernel-period-energy}
\end{equation}
The integrand is $2\pi n$-periodic. Since $p$ is an integer,
integration over any interval of length $2p\pi n$ therefore gives
at most $2\pi p$.

The support of $L_{n,s,\boldsymbol\ell}$ is contained in
$[-\pi n,\pi n]^p$, so $|y|\le p\pi n$ on that support.
Then, with further explanations given below, 
\begin{align}
\|L_{n,s,\boldsymbol\ell}\mathbf1_{A_n^c}\|_{L^2(\mathbb R^p)}^2
&\le
M^{|\boldsymbol\ell|}
\int_{\substack{|x_a|\le\pi n,\ a=1,\ldots,p\\
                 |x_1+\cdots+x_p|>\pi n}}
|\widehat K_{n,s}(x_1+\cdots+x_p)|^2
\prod_{a=1}^p\kappa_{k_a}|x_a|^{-\alpha_{k_a}}\,dx
\\
 &\le M^{|\boldsymbol\ell|}
 \int_{\pi n<|y|\le p\pi n}|\widehat K_{n,s}(y)|^2
 \Psi_{\boldsymbol\ell}(y)\,dy\label{eq:kernel-other-periods-2}\\
 &\le M^{|\boldsymbol\ell|} \Psi_{\boldsymbol\ell}(1) (\pi n)^{1-2H(\boldsymbol\ell)}
       \int_{-p\pi n}^{p\pi n}|\widehat K_{n,s}(y)|^2\,dy\label{eq:kernel-other-periods-1}\\
 &\le\pi pM^{|\boldsymbol\ell|}B_{\boldsymbol\ell}n^{1-2H(\boldsymbol\ell)}
 \to0.
 \label{eq:kernel-other-periods}
\end{align}
For the first inequality, we used \eqref{eq:kernel-factorization+1} and Assumption~\ref{ass:trace-admissibility}\ref{ass:coordinatewise-lrd}. The inequality \eqref{eq:kernel-other-periods-2} is due to substituting $
x_p=y-\sum_{a=1}^{p-1}x_a$, and using the definition of $\Psi_{\boldsymbol\ell}$ in \eqref{eq:psi-form}.
For \eqref{eq:kernel-other-periods-1}, we use \eqref{eq:power-convolution} of Lemma \ref{lem:envelope-finiteness} and
$|y|>\pi n$. The last inequality follows from
\eqref{eq:kernel-period-energy},
$\Psi_{\boldsymbol\ell}(1)\le B_{\boldsymbol\ell}/2$, and
$\pi^{1-2H(\boldsymbol\ell)}\le1$. The bound is uniform in $s\in[0,1]$ and
tends to zero because $2H(\boldsymbol\ell)-1>0$.
Combining \eqref{eq:kernel-central-convergence} and
\eqref{eq:kernel-other-periods}, the triangle inequality gives
\begin{equation}
\label{eq:scalar-kernel-convergence}
L_{n,s,\boldsymbol\ell}
\to
\kappa^{\boldsymbol\ell/2}K_{s,\boldsymbol\ell}
 \quad\text{in }L^2(\mathbb R^p).
\end{equation}

\textit{Proof of \eqref{eq:uniform-kernel-norm}:} By the isometry~\eqref{eq:labelled-isometry}, the coordinate spectral representation under
Assumption~\ref{ass:trace-admissibility}\ref{ass:diagonal-density} and the
change of variables $x_a=n\omega_a$,
\begin{equation}
\label{eq:kernel-variance-identity}
\|L_{n,s,\boldsymbol\ell}\|_{L^2(\mathbb R^p)}^2
=
\frac{n^{-2H(\boldsymbol\ell)}}{\boldsymbol\ell!}   \Var\left(\sum_{t=1}^{\lfloor ns\rfloor}H_{\boldsymbol\ell}[X_t]\right)
 \le\frac{C}{2 H(\boldsymbol\ell) -1}\left(\frac{\lfloor ns\rfloor}{n}\right)^{2H(\boldsymbol\ell)}
 \le\frac{C}{2H(\boldsymbol\ell)-1}.
\end{equation}
Here~\eqref{eq:spectral-block-bound} is applied with $d=H(\boldsymbol\ell)$. Fatou's lemma and the convergence in \eqref{eq:scalar-kernel-convergence} give
\[
\|\kappa^{\boldsymbol\ell/2}K_{s,\boldsymbol\ell}\|_{L^2(\mathbb R^p)}^2
 \le\liminf_{n\to\infty}\|L_{n,s,\boldsymbol\ell}\|_{L^2(\mathbb R^p)}^2
 \le\frac{C}{2H(\boldsymbol\ell)-1}.
\]
Both estimates hold for every $s\in[0,1]$, proving
\eqref{eq:uniform-kernel-norm}. The variance lemma was proved
directly from the spectral densities and does not use kernel
convergence or any coefficient summability.
\end{proof}

\subsection{Proof of Theorem~\ref{thm:operator-normalized-full}}
\label{sec:proof-main}
\begin{proof}
We combine
Propositions 
\ref{prop:operator-reduction-principle} and \ref{prop:operator-kernel-convergence}. We first prove convergence of finite-dimensional distributions in a product space of $\mathbb H_2$.
We then consider the continuous
interpolations \eqref{eq:continuous-interpolation} in $C([0,1];\mathbb{H}_2)$ and tightness.

\medskip\noindent\textit{Convergence of finite-dimensional distributions:}
Recall $G_{\mathcal D}$ from \eqref{eq:dominant-block}, and write
\begin{align}
 A_n(s)=n^{-D}\sum_{t=1}^{\lfloor ns\rfloor}G_{\mathcal D}[X_t],
\qquad
 R_n(s)=n^{-D}\sum_{t=1}^{\lfloor ns\rfloor}
(G[X_t]-G_{\mathcal D}[X_t]).
 \label{eq:main-decomposition}
\end{align}
It then suffices to prove convergence in distribution of $A_n(s)$ and that $R_n(s)$ is negligible to infer convergence of $A_n+R_n$ matching the prelimit in
\eqref{eq:operator-normalized-convergence}.

For $N\ge1$, define the finite projections, using the notation \eqref{eq:def:projections:G} with $\mathcal{D} = (\mathcal{D}_i)_{i \geq 1}$,
\begin{align*}
 A_n^{(N)}(s)
 =n^{-D}\sum_{t=1}^{\lfloor ns\rfloor}G_{\mathcal D}^{(N)}[X_t],
 \qquad
 Z_D^{(N)}(s)
=\sum_{(i,\boldsymbol\ell)\in\mathcal F_N(\mathcal D)}
 c_{i,\boldsymbol\ell}\kappa^{\boldsymbol\ell/2}
 I_{|\boldsymbol\ell|}(J_{\boldsymbol\ell}K_{s,\boldsymbol\ell})v_i.
\end{align*}
We first prove convergence for fixed $N$.
For every coordinate spectral measure $\widehat W_j$, define
\[
 \widehat W_{j,n}(B)=n^{1/2}\widehat W_j(B/n),
 \qquad B\subset\mathbb R\text{ Borel with finite Lebesgue measure}.
\]
The rescaled measures have Lebesgue control measure, since
$n|B/n|=|B|$, and retain Hermitian symmetry and independence
across coordinate labels. Hence the finitely many rescaled measures
have the same joint law as the original measures.
Applying the change of variables $x_a=n\omega_a$ simultaneously to
the finite collection of multiple integrals gives
\begin{equation}
\label{eq:main-coupled-dominant}
 A_n^{(N)}(s)
 \ \stackrel{\operatorname{f.d.d.}}{=}\
 \widehat A_n^{(N)}(s)
 \doteq\sum_{(i,\boldsymbol\ell)\in\mathcal F_N(\mathcal D)}
 c_{i,\boldsymbol\ell}
 I_{|\boldsymbol\ell|}(J_{\boldsymbol\ell}L_{n,s,\boldsymbol\ell})v_i.
\end{equation}
Here we use \eqref{eq:hermite-spectral-integral} and
$H(\boldsymbol\ell)=d_i$ on $\mathcal D_i$.
This is the usual joint change-of-measure property of multiple
Wiener integrals; see \cite[Appendix~B.2, equation~(B.2.15)]{PipirasTaqqu17}.
The same rescaling is used for all selected times, labels and chaos
orders, so \eqref{eq:main-coupled-dominant} is a joint equality in law.

Fix $s_1,\ldots,s_m\in[0,1]$. The labeled isometry and
orthogonality give
\begin{align}
 \E\bigl\|(
 \widehat A_n^{(N)}(s_a)-Z_D^{(N)}(s_a))_{a=1}^m
 \bigr\|_{\mathbb H_2^m}^2
 =
\sum_{a=1}^m\sum_{(i,\boldsymbol\ell)\in\mathcal F_N(\mathcal D)}
\boldsymbol\ell!|c_{i,\boldsymbol\ell}|^2
 \|L_{n,s_a,\boldsymbol\ell}
      -\kappa^{\boldsymbol\ell/2}K_{s_a,\boldsymbol\ell}\|_{L^2}^2
 \to0.
 \label{eq:main-coupled-error}
\end{align}
There are only finitely many summands, and each tends to zero by
\eqref{eq:scalar-kernel-convergence}. Thus, for fixed $N$,
$(A_n^{(N)}(s_a))_{a=1}^m$ converges in distribution to
$(Z_D^{(N)}(s_a))_{a=1}^m$ in the finite-dimensional subspace
$\operatorname{span}\{v_1,\ldots,v_N\}^m$.

It remains to remove the projection. By
Lemma~\ref{lem:finite-hermite-projection} and
\eqref{eq:uniform-kernel-norm}, respectively,
\begin{align}
 \sup_{n\ge1}\sup_{s\in[0,1]}
 \E\|A_n(s)-A_n^{(N)}(s)\|_{\mathbb H_2}^2
 &\le C \sum_{\substack{i\ge1,\ \boldsymbol\ell\in\mathcal D_i\\
                   (i,\boldsymbol\ell)\notin\mathcal F_N(\mathcal D)}}
 \frac{\boldsymbol\ell!|c_{i,\boldsymbol\ell}|^2}{2d_i-1},\notag\\
 \sup_{s\in[0,1]}
 \E\|Z_D(s)-Z_D^{(N)}(s)\|_{\mathbb H_2}^2
 &\le C \sum_{\substack{i\ge1,\ \boldsymbol\ell\in\mathcal D_i\\
                   (i,\boldsymbol\ell)\notin\mathcal F_N(\mathcal D)}}
 \frac{\boldsymbol\ell!|c_{i,\boldsymbol\ell}|^2}{2d_i-1}.
 \label{eq:main-projection-errors}
\end{align}
For the second inequality, use the labeled isometry and
$H(\boldsymbol\ell)=d_i$ on $\mathcal D_i$.
By Parseval's identity and Assumption~\ref{ass:bigAssump2}\ref{ass:summability},
the full weighted coefficient sum is finite. Since the sets
$\mathcal F_N(\mathcal D)$ increase to all dominant indices,
dominated convergence shows that both right-hand sides tend to zero
as $N\to\infty$.
Summing over the $m$ times and applying Markov's inequality verifies
the approximation condition of \cite[Theorem~3.2]{billingsley1999convergence}.
That theorem therefore yields
\begin{equation}
 \label{eq:main-product-convergence}
 (A_n(s_a))_{a=1}^m\xrightarrow{d}(Z_D(s_a))_{a=1}^m
 \quad\text{in }\mathbb H_2^m.
\end{equation}

We turn to the remainder $R_n(s)$ in \eqref{eq:main-decomposition}. Proposition~\ref{prop:operator-reduction-principle}
implies, for every $\varepsilon>0$,
\begin{align}
\Prob\bigl(\|(R_n(s_a))_{a=1}^m\|_{\mathbb H_2^m}>\varepsilon\bigr) \le\varepsilon^{-2}\sum_{a=1}^m\E\|R_n(s_a)\|_{\mathbb H_2}^2
\le\varepsilon^{-2}m\sup_{s\in[0,1]}\E\|R_n(s)\|_{\mathbb H_2}^2
\to0,
\label{eq:main-remainder-vector}
\end{align}
where the first inequality is Markov's inequality. Combining
\eqref{eq:main-decomposition}, \eqref{eq:main-product-convergence},
and \eqref{eq:main-remainder-vector}, Slutsky's lemma in
$\mathbb H_2^m$ proves \eqref{eq:operator-normalized-convergence}.

\medskip\noindent\textit{Continuous interpolations:}
Recall \eqref{eq:continuous-interpolation}, with $T_n=n^{-D}$.
Since Assumption~\ref{ass:bigAssump2}\ref{ass:operator-standing} gives $d_i>1/2$ for every $i$, including the prescribed values on
zero output directions, $\|n^{-D}\|_{\mathcal L(\mathbb H_2)}\le n^{-1/2}$.
Stationarity therefore gives, for every fixed $s\in[0,1]$,
\begin{align}
\E\left\|\widetilde S_n(s)
-n^{-D}\sum_{t=1}^{\lfloor ns\rfloor}G[X_t]\right\|_{\mathbb H_2}^2
\le\|n^{-D}\|_{\mathcal L(\mathbb H_2)}^2\E\|G[X_0]\|_{\mathbb H_2}^2
\le n^{-1}\E\|G[X_0]\|_{\mathbb H_2}^2\to0.
 \label{eq:main-interpolation-error}
\end{align}
Summing this estimate over $s_1,\ldots,s_m$ shows that the
interpolations have the same limiting finite-dimensional distributions.

\medskip\noindent\textit{Tightness:}
We apply the tightness criterion in
\cite[Proposition~4.1]{KraetschmerUrusov2023} with parameter space
$[0,1]$ and state space $\mathbb H_2$.
The marginal laws are tight at every fixed time by the convergence
already proved. It remains to establish a uniform increment bound
with a power strictly greater than one.

We impose the additional condition $d_*>1/2$ in
Theorem~\ref{thm:operator-normalized-full}.
For every integer $a$ and every integer $m\ge1$, stationarity
and the finite-projection argument in the proof of
Lemma~\ref{lem:finite-hermite-projection} give
\begin{align}
\E\left|\sum_{t=a+1}^{a+m}
\langle G[X_t],v_i\rangle_{\mathbb H_2}\right|^2
&=\sum_{p=1}^{\infty}\sum_{\boldsymbol\ell\in\mathcal L_p}
|c_{i,\boldsymbol\ell}|^2
\Var\left(\sum_{t=1}^m H_{\boldsymbol\ell}[X_t]\right)
\notag\\
&\le \frac{C}{2d_i-1}
\E|\langle G[X_0],v_i\rangle_{\mathbb H_2}|^2\,m^{2d_i}.
\label{eq:tightness-block-estimate}
\end{align}
For the inequality, every index with $c_{i,\boldsymbol\ell}\ne0$
satisfies $H(\boldsymbol\ell)\le d_i$ by the definition of $d_i$
and Assumption~\ref{ass:bigAssump2}\ref{ass:operator-standing}. We therefore apply
\eqref{eq:spectral-block-bound} of
Lemma~\ref{lem:spectral-variance}, and then Parseval's
identity for the scalar Hermite expansion of $\langle G,v_i\rangle$.

Let $U_{n,i}(s)=\langle\widetilde S_n(s),v_i\rangle_{\mathbb H_2}$, with
$T_n=n^{-D}$, and put $h=t-s>0$. We may take the constant $C$ in
\eqref{eq:tightness-block-estimate} to be at least one.
For grid times $s,t\in\{0,1/n,\ldots,1\}$ with $s<t$, we have
\begin{align}
\|U_{n,i}(t)-U_{n,i}(s)\|_{L^2(\Omega)}
&=
n^{-d_i}  \left\|\sum_{r=ns+1}^{nt}\langle G[X_r],v_i\rangle_{\mathbb H_2}\right\|_{L^2(\Omega)}
\\&\le
 \left(\frac{C}{2d_i-1}\right)^{1/2}
 \|\langle G[X_0],v_i\rangle_{\mathbb H_2}\|_{L^2(\Omega)}\,h^{d_i},
\end{align}
where the equality follows from the definition of the interpolation at
grid times. The inequality is
\eqref{eq:tightness-block-estimate} with block length $n(t-s)$.
For arbitrary $s,t$ with $h\ge1/n$, put
$j=\lfloor ns\rfloor$ and $k=\lfloor nt\rfloor$, so that $k\ge j+1$.
We split the increment as follows
\begin{align}
 &
 \|U_{n,i}(t)-U_{n,i}(s)\|_{L^2(\Omega)}
 \\&\le
 \|U_{n,i}(t)-U_{n,i}(k/n)\|_{L^2(\Omega)}
+
 \|U_{n,i}(k/n)-U_{n,i}((j+1)/n)\|_{L^2(\Omega)}
\\&\hspace{6cm}+
 \|U_{n,i}((j+1)/n)-U_{n,i}(s)\|_{L^2(\Omega)}
 \notag\\
 &\le
 \left(\frac{C}{2d_i-1}\right)^{1/2}
 \|\langle G[X_0],v_i\rangle_{\mathbb H_2}\|_{L^2(\Omega)}
 \left[
 \frac{nt-k}{n^{d_i}}
 +\left(\frac{k-j-1}{n}\right)^{d_i}
 +\frac{j+1-ns}{n^{d_i}}
 \right]
 \notag\\
 &\le
 3\left(\frac{C}{2d_i-1}\right)^{1/2}
 \|\langle G[X_0],v_i\rangle_{\mathbb H_2}\|_{L^2(\Omega)}\,h^{d_i}.
\end{align}
The first inequality is the triangle inequality. For the second,
linear interpolation and the case $m=1$ of
\eqref{eq:tightness-block-estimate} control the two partial
increments, while the grid-time estimate controls the middle block.
For the last inequality, each partial-increment coefficient is at
most one, $n^{-d_i}\le h^{d_i}$, and
$(k-j-1)/n\le h$.

If $0<h<1/n$, we use the piecewise constant slope of the linear
interpolation to get, with further explanations below,
\begin{align}
 \|U_{n,i}(t)-U_{n,i}(s)\|_{L^2(\Omega)}
 &=n^{1-d_i}
 \left\|\int_s^t
 \langle G[X_{\lfloor nv\rfloor+1}],v_i\rangle_{\mathbb H_2}\,dv
 \right\|_{L^2(\Omega)}
 \label{eq:tightness-increment-estimate-3-3}\\
 &\le n^{1-d_i}\int_s^t
 \|\langle G[X_{\lfloor nv\rfloor+1}],v_i\rangle_{\mathbb H_2}\|_{L^2(\Omega)}\,dv
 \label{eq:tightness-increment-estimate-3-2}\\
 &=n^{1-d_i}h\,
 \|\langle G[X_0],v_i\rangle_{\mathbb H_2}\|_{L^2(\Omega)}
 \label{eq:tightness-increment-estimate-3-1}\\
 &\le h^{d_i}\|\langle G[X_0],v_i\rangle_{\mathbb H_2}\|_{L^2(\Omega)},
 \label{eq:tightness-increment-estimate-3}
\end{align}
where \eqref{eq:tightness-increment-estimate-3-3} follows by integrating the piecewise constant slope of the linear interpolation defined in \eqref{eq:continuous-interpolation}, \eqref{eq:tightness-increment-estimate-3-2} uses Minkowski's integral inequality, and \eqref{eq:tightness-increment-estimate-3-1} uses stationarity. The final inequality \eqref{eq:tightness-increment-estimate-3} factors out $h^{d_i}$,
and uses $nh<1$ and $d_i<1$, the latter by
Assumption~\ref{ass:bigAssump2}\ref{ass:operator-standing}.
Since $C\ge1$ and $2d_i-1<1$, both cases are bounded by the estimate
with factor $3(C/(2d_i-1))^{1/2}$ above. Consequently,
\begin{align}
 \E\|\widetilde S_n(t)-\widetilde S_n(s)\|_{\mathbb H_2}^2
 &=\sum_{i=1}^\infty \|U_{n,i}(t)-U_{n,i}(s)\|_{L^2(\Omega)}^2
 \notag\\
 &\le9C\sum_{i=1}^\infty
 \frac{\E|\langle G[X_0],v_i\rangle_{\mathbb H_2}|^2}{2d_i-1}
 |t-s|^{2d_i}
 \label{eq:tightness-increment-estimate-1}\\
 &\le9C|t-s|^{2d_*}
 \sum_{i=1}^\infty \frac{\E|\langle G[X_0],v_i\rangle_{\mathbb H_2}|^2}{2d_i-1}.
 \label{eq:tightness-increment-estimate}
\end{align}
The equality follows from Parseval's identity and Tonelli's theorem.
The inequality \eqref{eq:tightness-increment-estimate-1} uses the preceding two cases,
and \eqref{eq:tightness-increment-estimate} uses $|t-s|\le1$ and $d_i\ge d_*$ on active
coordinates. Inactive coordinates contribute zero. The coefficient
sum is finite by Assumption \ref{ass:bigAssump2}\ref{ass:summability}.

Since $2d_*>1$, \cite[Proposition~4.1]{KraetschmerUrusov2023}
gives tightness in $C([0,1];\mathbb H_2)$.
Its covering-number condition holds on $[0,1]$ with exponent one;
the increment moment has order two and time exponent $2d_*$.
The state space is complete and separable, and the interpolations
have continuous paths, so all its hypotheses are satisfied.

For completeness, the same increment estimate also gives a
continuous version of the limiting process. By the finite-dimensional
convergence and the Portmanteau theorem,
\[
 \E\|Z_D(t)-Z_D(s)\|_{\mathbb H_2}^2
 \le 9C|t-s|^{2d_*}
 \sum_{i=1}^{\infty}
 \frac{\E|\langle G[X_0],v_i\rangle_{\mathbb H_2}|^2}{2d_i-1}.
\]
Hence \cite[Theorem~1.1]{KraetschmerUrusov2023} supplies a continuous
modification. Finally, tightness and the finite-dimensional
convergence, together with \cite[Lemma~16.2]{Kallenberg2002}, yield
\[
 \widetilde S_n \xrightarrow{d} Z_D
 \quad\text{in }C([0,1];\mathbb H_2),
\]
where $Z_D$ denotes this continuous version.
\end{proof}

\subsection{Proof of Corollary~\ref{thm:common-hermite-limit}}
\label{sec:proof-common}
\begin{proof}
Recall the expansion $G[x] = \sum_{i\ge1}\sum_{\ell\in\mathcal L}
c_{i,\ell}H_\ell[x]\otimes v_i$ from \eqref{eq:full-hermite-expansion}.
We aim to apply Theorem~\ref{thm:operator-normalized-full} after expanding $G[x]$ with respect to a basis that is possibly different than $\{v_i\}_{i \geq 1}$ and in which every coordinate of $G$ has a nonzero term of order $q$. We then verify the assumption of Theorem~\ref{thm:operator-normalized-full}. For every $\boldsymbol\ell\in\mathcal L$, denote by
\[
  c_{\boldsymbol\ell}
  \doteq\sum_{j=1}^{\infty}c_{j,\boldsymbol\ell}v_j
  \in\mathbb H_2
\]
the coefficient vector in \eqref{eq:full-hermite-expansion}. Since $G$ has
Hermite rank $q$, there exists
$\boldsymbol\ell_0\in\mathcal L_q$ such that
$c_{\boldsymbol\ell_0}\ne0$.
To construct the required basis, choose
$\{a_i\}_{i\ge1}\in\ell^2$ with $a_i>0$ for every $i\ge1$ and
$\sum_i a_i^2=1$, and put $u=\sum_i a_i v_i$. By extending $u$
and $c_{\boldsymbol\ell_0}/\|c_{\boldsymbol\ell_0}\|_{\mathbb H_2}$
to orthonormal bases, we obtain a surjective linear isometry $U$
of $\mathbb H_2$ satisfying
$Uu=c_{\boldsymbol\ell_0}/\|c_{\boldsymbol\ell_0}\|_{\mathbb H_2}$.
Let $w_i=Uv_i$, and denote by $\widetilde c_{i,\boldsymbol\ell}$
the Hermite coefficients in the orthonormal basis $\{w_i\}_{i\ge1}$.
For every $i\ge1$ and $\boldsymbol\ell\in\mathcal L$, we have
\begin{equation}
\begin{gathered}
\label{eq:defofciltildebasis}
  \widetilde c_{i,\boldsymbol\ell}
  =\langle c_{\boldsymbol\ell},w_i\rangle_{\mathbb H_2}
   =\sum_{j=1}^{\infty}
     c_{j,\boldsymbol\ell}\langle v_j,w_i\rangle_{\mathbb H_2},
  \quad
  \widetilde c_{i,\boldsymbol\ell_0}
  =\|c_{\boldsymbol\ell_0}\|_{\mathbb H_2}
    \langle Uu,Uv_i\rangle_{\mathbb H_2}
   =\|c_{\boldsymbol\ell_0}\|_{\mathbb H_2}a_i>0,
  \\
  \sum_{i=1}^{\infty}|\widetilde c_{i,\boldsymbol\ell}|^2
  =\|c_{\boldsymbol\ell}\|_{\mathbb H_2}^2  =\sum_{j=1}^{\infty}|c_{j,\boldsymbol\ell}|^2,
\end{gathered}
\end{equation}
where the representation of $\widetilde c_{i,\boldsymbol\ell_0}$ uses the isometry of $U$, and the last line
follows from Parseval's identity. In this basis,
\[
\begin{aligned}
G[X_t] =\sum_{i=1}^{\infty}
    \langle G[X_t],w_i\rangle_{\mathbb H_2}\,w_i
=\sum_{i=1}^{\infty}
    \left(
      \sum_{\boldsymbol\ell\in\mathcal L}
      \widetilde c_{i,\boldsymbol\ell}H_{\boldsymbol\ell}[X_t]
    \right)w_i
=\sum_{i=1}^{\infty}\sum_{\boldsymbol\ell\in\mathcal L}
    \widetilde c_{i,\boldsymbol\ell}
    H_{\boldsymbol\ell}[X_t]\,w_i,
\end{aligned}
\]
with convergence in $L^2(\Omega;\mathbb H_2)$. Thus the output
basis and coefficients are replaced by $w_i$ and
$\widetilde c_{i,\boldsymbol\ell}$, respectively. Moreover, the coefficient norm identity above
shows that a multi-index is active in the new basis if and only
if it is active in the original basis. Consequently, the active
input coordinates are unchanged, and
Assumption~\ref{ass:trace-admissibility} remains valid.

To verify Assumption~\ref{ass:bigAssump2} in the basis $\{w_i\}_{i \geq 1}$, denote by $\widetilde d_i$, $\widetilde{\mathcal D}_i$, and
$\widetilde D$ the corresponding dominant exponents, dominant
sets, and operator, respectively, as defined in
\eqref{eq:dominant-exponent-and-set} and
\eqref{eq:operator-normalization-D}. By
Assumption~\ref{ass:common-active-exponent}\ref{ass:common-active-dom-chaos},
all input coordinates in the support of an active index of
order $q$ have exponent $\alpha$. Hence, for every
$\boldsymbol\ell\in\mathcal L_q$ with $c_{\boldsymbol\ell}\ne0$,
\begin{align*}
\beta(\boldsymbol\ell)
=\sum_j\ell_j(1-\alpha_j)
=(1-\alpha)|\boldsymbol\ell|
=q(1-\alpha),
\quad
H(\boldsymbol\ell)
=1-\frac q2(1-\alpha)
=H\in(1/2,1),
\end{align*}
and $K_{s,\boldsymbol\ell} = K_s$ by \eqref{eq:operator-limit-kernel}.
Moreover,
Assumption~\ref{ass:common-active-exponent}\ref{ass:common-active-high-chaos}
gives $H(\boldsymbol\ell)<H$ for every active higher-order index.
Since $\widetilde c_{i,\boldsymbol\ell_0}>0$ for every $i\ge1$,
the dominant exponent in each output coordinate is attained at
$\boldsymbol\ell_0$. Consequently,
\[
  \widetilde d_i=H,
  \qquad
  \widetilde{\mathcal D}_i
  =\{\boldsymbol\ell\in\mathcal L_q:
       \widetilde c_{i,\boldsymbol\ell}\ne0\},
  \qquad
  \widetilde D=H I_{\mathbb H_2},
\]
which verifies
Assumption~\ref{ass:bigAssump2}\ref{ass:operator-standing}.
Furthermore, Parseval's identity and
Assumption~\ref{ass:common-active-exponent}\ref{ass:GinL2} yield
\begin{align}
  \sum_{i=1}^{\infty}
  \frac{\E|\langle G[X_0],w_i\rangle_{\mathbb H_2}|^2}
       {2\widetilde d_i-1}
  &=\frac{1}{2H-1}
    \sum_{i=1}^{\infty}\sum_{\boldsymbol\ell\in\mathcal L}
    \boldsymbol\ell!|\widetilde c_{i,\boldsymbol\ell}|^2
    \notag\\
  &=\frac{1}{2H-1}
    \sum_{\boldsymbol\ell\in\mathcal L}
    \boldsymbol\ell!\|c_{\boldsymbol\ell}\|_{\mathbb H_2}^2
 =\frac{\E\|G[X_0]\|_{\mathbb H_2}^2}{2H-1}<\infty,
  \label{eq:common-summability-verification}
\end{align}
which verifies
Assumption~\ref{ass:bigAssump2}\ref{ass:summability}.

Since $\widetilde D=H I_{\mathbb H_2}$ and
$\inf_i\widetilde d_i=H>1/2$,
Theorem~\ref{thm:operator-normalized-full}, applied in the basis
$\{w_i\}_{i \geq 1}$, gives both finite-dimensional convergence and weak
convergence of the continuous interpolations with
$T_n=n^{-H} I_{\mathbb H_2}$. Denote by $\widetilde Z$ the resulting limit.
It remains to identify $\widetilde Z$ with $Z_H^{(q)}$. For every
$s\in[0,1]$, \eqref{eq:operator-limit-process} and the kernel
identity above give
\begin{align*}
\widetilde Z(s) 
&=
\sum_{i=1}^{\infty}\sum_{\boldsymbol\ell\in\mathcal L_q}
\widetilde c_{i,\boldsymbol\ell}\kappa^{\boldsymbol\ell/2}
I_q(J_{\boldsymbol\ell}K_s)w_i
=
\sum_{\boldsymbol\ell\in\mathcal L_q} \kappa^{\boldsymbol\ell/2}I_q(J_{\boldsymbol\ell}K_s)
    \left(
      \sum_{i=1}^{\infty}\widetilde c_{i,\boldsymbol\ell}w_i
    \right)
\\
&=\sum_{\boldsymbol\ell\in\mathcal L_q}
    \kappa^{\boldsymbol\ell/2}I_q(J_{\boldsymbol\ell}K_s)
    \left(
      \sum_{i=1}^{\infty}
      \langle c_{\boldsymbol\ell},w_i\rangle_{\mathbb H_2}w_i
    \right)
=
\sum_{\boldsymbol\ell\in\mathcal L_q}
    \kappa^{\boldsymbol\ell/2}I_q(J_{\boldsymbol\ell}K_s)
    c_{\boldsymbol\ell}
  \\
  &=\sum_{\boldsymbol\ell\in\mathcal L_q}
    \kappa^{\boldsymbol\ell/2}I_q(J_{\boldsymbol\ell}K_s)
    \left(
      \sum_{j=1}^{\infty}c_{j,\boldsymbol\ell}v_j
    \right)
=\sum_{j=1}^{\infty}\sum_{\boldsymbol\ell\in\mathcal L_q}
    c_{j,\boldsymbol\ell}\kappa^{\boldsymbol\ell/2}
    I_q(J_{\boldsymbol\ell}K_s)v_j
=Z_H^{(q)}(s),
\end{align*}
where we used \eqref{eq:defofciltildebasis}.
This identification yields \eqref{eq:common-hermite-convergence}
and weak convergence in $C([0,1];\mathbb H_2)$ to a continuous
version of $Z_H^{(q)}$.
\end{proof}

\section{Additional results and their proofs} \label{sec:auxiliary-lemmas}

The following lemma provides a spectral comparison for Hermite
partial sums.

\begin{lemma}
\label{lem:hermite-variance-comparison}
Suppose  Assumptions~\ref{ass:trace-admissibility}\ref{ass:diagonal-density}--\ref{ass:uniform-spectral-scale}.
There exists a constant $C\ge1$, depending only on $M$,
$\kappa_*$ and $\sup_j\alpha_j$, such that, for every
$\boldsymbol\ell\in\mathcal L_p$, $p\ge1$, supported on active
coordinates, and all integers $m,r\ge1$,
\begin{equation}
\label{eq:hermite-partial-comparison}
\Var\left(\sum_{t=1}^m H_{\boldsymbol\ell}[X_t]\right)
\le
\boldsymbol\ell! C^{p\wedge r}
\sum_{|h|<m}(m-|h|)
(1+|h|)^{-\sum_{a=1}^{p\wedge r}(1-\alpha_{k_a})},
\end{equation}
where $\mathbf k=\mathbf k(\boldsymbol\ell)$.
\end{lemma}

\begin{proof}
By Assumption~\ref{ass:trace-admissibility}\ref{ass:coordinatewise-lrd}, we have 
$\widetilde f_j(\omega)\le
g_j(\omega)\doteq M \kappa_j|\omega|^{-\alpha_j}$
for every active coordinate $j$.
For $u\in L^1([-\pi,\pi])$, write
\[
\widehat u(h)=\int_{-\pi}^{\pi}e^{ih\omega}u(\omega)\,d\omega.
\]
We first note that, uniformly over active coordinates, 
\begin{equation}
\label{eq:hermite-comparison-fourier}
|\widehat g_j(h)|
\le C(1+|h|)^{-(1-\alpha_j)},
\qquad
|\widehat{\widetilde f}_j(h)|\le1,
\qquad h\in\mathbb Z.
\end{equation}
Indeed, direct integration gives
$\widehat g_j(0)=
2M\kappa_j\pi^{1-\alpha_j}/(1-\alpha_j)$.
For $h\ne0$, evenness and a change of variables give
\[
|\widehat g_j(h)|
=
2M\kappa_j|h|^{\alpha_j-1}
\left|\int_0^{\pi|h|}u^{-\alpha_j}\cos u\,du\right|
\le
2M\kappa_j|h|^{\alpha_j-1}
\left(\frac{1}{1-\alpha_j}+2\right).
\]
The last inequality follows by splitting the integral at $1$
and integrating by parts on $[1,\pi|h|]$.
The bounds $\kappa_j\le\kappa_*$ from Assumption~\ref{ass:trace-admissibility}\ref{ass:uniform-spectral-scale} and
$\sup_j\alpha_j<1$ from Assumption~\ref{ass:trace-admissibility}\ref{ass:coordinatewise-lrd} yield the first estimate in
\eqref{eq:hermite-comparison-fourier}.
The second follows from
$\int_{-\pi}^{\pi}\widetilde f_j(\omega)\,d\omega=1$.

Put $q=p\wedge r$ and define $D_m(y)=\sum_{t=1}^m e^{ity}$.
Expanding the variance and substituting the Fourier
representations yields
\begin{align*}
\Var\left(\sum_{t=1}^m H_{\boldsymbol\ell}[X_t]\right)
&=
\boldsymbol\ell!
\int_{[-\pi,\pi]^p}
|D_m(\omega_1+\cdots+\omega_p)|^2
\prod_{a=1}^p\widetilde f_{k_a}(\omega_a)
\,d\boldsymbol\omega\\
&\le
\boldsymbol\ell!
\int_{[-\pi,\pi]^p}
|D_m(\omega_1+\cdots+\omega_p)|^2
\prod_{a=1}^q g_{k_a}(\omega_a)
\prod_{a=q+1}^p\widetilde f_{k_a}(\omega_a)
\,d\boldsymbol\omega\\
&\le
\boldsymbol\ell!
\sum_{|h|<m}(m-|h|)
\prod_{a=1}^q|\widehat g_{k_a}(h)|
\prod_{a=q+1}^p|\widehat{\widetilde f}_{k_a}(h)|.
\end{align*}
The first inequality follows from the nonnegativity of the
integrand. For the second, we expand $|D_m|^2$, integrate
each frequency separately, and take absolute values;
there are $m-|h|$ pairs of time indices with difference $h$.
All interchanges involve finite sums and integrable functions,
and empty products are understood to equal one.
Applying \eqref{eq:hermite-comparison-fourier} proves
\eqref{eq:hermite-partial-comparison}.
\end{proof}

The following lemma employs the bounds of Lemma \ref{lem:hermite-variance-comparison} to obtain a variance bound for Hermite partial sums that is uniform in the Hermite order.


\begin{lemma}
\label{lem:spectral-variance}
Suppose Assumptions~\ref{ass:trace-admissibility}\ref{ass:diagonal-density}--\ref{ass:uniform-spectral-scale}.
There exists a constant $C\ge1$, depending only on $M$, $\kappa_*$
and $\sup_j\alpha_j$, such that, for every
$\boldsymbol\ell\in\mathcal L_p$, $p\ge1$, supported on active
coordinates, every integer $m \geq 1$, and every 
$d\in(1/2,1)$ satisfying $H(\boldsymbol\ell)\le d$,
\begin{equation}
\label{eq:spectral-block-bound}
\Var\left(\sum_{t=1}^mH_{\boldsymbol\ell}[X_t]\right)
 \le\frac{C\boldsymbol\ell!}{2d-1}\,m^{2d}.
\end{equation}
\end{lemma}

\begin{proof}
Under Assumptions~\ref{ass:trace-admissibility}\ref{ass:diagonal-density}--\ref{ass:uniform-spectral-scale}, let $C_0\ge1$ be the constant in
Lemma~\ref{lem:hermite-variance-comparison}.
Fix $\boldsymbol\ell\in\mathcal L_p$ supported on active coordinates,
and write $\mathbf k=\mathbf k(\boldsymbol\ell)$ and
$\beta=\beta(\boldsymbol\ell)$.
Applying that lemma with $r=p$ gives the lag-sum estimate used in
\cite[proof of Lemma~3.1]{DehlingTaqqu89}:
\begin{align}
\Var\left(\sum_{t=1}^mH_{\boldsymbol\ell}[X_t]\right)
&\le \boldsymbol\ell!C_0^p m
\left(1+2\sum_{h=1}^{m-1}(1+h)^{-\beta}\right)
\notag\\
&\le \boldsymbol\ell!C_0^p m
\left(1+2\int_1^m x^{-\beta}\,dx\right).
\label{eq:hermite-sum-variance-1}
\end{align}
The first inequality separates the zero lag and uses $m-h\le m$.
The second follows because $x^{-\beta}$ is decreasing.

For uniformity in the chaos order, we use fixed-order domination
as in \cite[equation~(37)]{bai2019sensitivity}.
Put $\varepsilon=1-\sup_j\alpha_j>0$ by
Assumption~\ref{ass:trace-admissibility}\ref{ass:coordinatewise-lrd}, and choose an integer
$r_0$ depending only on $\varepsilon$ such that $r_0\varepsilon>1$.
For $p\ge r_0$, Lemma~\ref{lem:hermite-variance-comparison}
with $r=r_0$ gives
\begin{align}
\Var\left(\sum_{t=1}^mH_{\boldsymbol\ell}[X_t]\right)
&\le \boldsymbol\ell!C_0^{r_0}m
\left(1+2\sum_{h=1}^{\infty}(1+h)^{-r_0\varepsilon}\right)
\notag\\
&\le \boldsymbol\ell!C_0^{r_0}m
\left(1+\frac{2}{r_0\varepsilon-1}\right)
\doteq C_1\boldsymbol\ell!m.
\label{eq:high-chaos-uniform-bound}
\end{align}
Here we used
$\sum_{a=1}^{r_0}(1-\alpha_{k_a})\ge r_0\varepsilon$
and the integral comparison on $[1,\infty)$.
The constant $C_1$ depends only on $M$, $\kappa_*$ and
$\sup_j\alpha_j$.

Finally, let $d\in(1/2,1)$ satisfy $H(\boldsymbol\ell)\le d$.
For $p<r_0$, we have $\beta=2-2H(\boldsymbol\ell)\ge2-2d$,
so the first bound in \eqref{eq:hermite-sum-variance-1} yields
\begin{align*}
\Var\left(\sum_{t=1}^mH_{\boldsymbol\ell}[X_t]\right)
&\le \boldsymbol\ell!C_0^{r_0}m
\left(1+2\int_1^m x^{2d-2}\,dx\right)\\
&=\boldsymbol\ell!C_0^{r_0}m
\left(1+2\frac{m^{2d-1}-1}{2d-1}\right)
\le \frac{3C_0^{r_0}\boldsymbol\ell!}{2d-1}m^{2d}.
\end{align*}
We used $C_0^p\le C_0^{r_0}$ and $0<2d-1<1$.
For $p\ge r_0$, \eqref{eq:high-chaos-uniform-bound} gives
\[
\Var\left(\sum_{t=1}^mH_{\boldsymbol\ell}[X_t]\right)
\le C_1\boldsymbol\ell!m
\le \frac{C_1\boldsymbol\ell!}{2d-1}m^{2d}.
\]
Choosing $C=\max\{C_0,C_1,3C_0^{r_0}\}$ proves 
\eqref{eq:spectral-block-bound}.
\end{proof}

The following lemma provides a uniform bound for the projection remainder.

\begin{lemma}
\label{lem:finite-hermite-projection}
Suppose Assumptions~\ref{ass:trace-admissibility}\ref{ass:diagonal-density}--\ref{ass:uniform-spectral-scale} and
$H(\boldsymbol\ell)\le d_i$ for every $i\ge1$ and
$\boldsymbol\ell\in\mathcal Q_i$.
Then there is a constant $C>0$, such that, for all $n,N\ge1$,
\begin{equation}
\label{eq:projection-remainder-bound}
\sup_{s\in[0,1]}\E\left\|
n^{-D}\sum_{t=1}^{\lfloor ns\rfloor}
(G_{\mathcal Q}-G_{\mathcal Q}^{(N)})[X_t]
\right\|_{\mathbb H_2}^2
\le 
C \sum_{\substack{i\ge1,\ \boldsymbol\ell\in\mathcal Q_i\\
(i,\boldsymbol\ell)\notin\mathcal F_N(\mathcal Q)}}
\frac{\boldsymbol\ell!|c_{i,\boldsymbol\ell}|^2}{2d_i-1}.
\end{equation}
\end{lemma}

\begin{proof}
For $L\ge N$ and $0\le m\le n$, Lemma \ref{lem:hermite-basis-detail}, gives
\begin{align*}
&\E\left\|n^{-D}\sum_{t=1}^m
(G_{\mathcal Q}^{(L)}-G_{\mathcal Q}^{(N)})[X_t]
\right\|_{\mathbb H_2}^2\\
&=
\sum_{(i,\boldsymbol\ell)\in
\mathcal F_L(\mathcal Q)\setminus\mathcal F_N(\mathcal Q)}
|c_{i,\boldsymbol\ell}|^2 n^{-2d_i}
\Var\left(\sum_{t=1}^mH_{\boldsymbol\ell}[X_t]\right)
 \le
C \sum_{\substack{i\ge1,\ \boldsymbol\ell\in\mathcal Q_i\\
(i,\boldsymbol\ell)\notin\mathcal F_N(\mathcal Q)}}
\frac{\boldsymbol\ell!|c_{i,\boldsymbol\ell}|^2}{2d_i-1}.
\end{align*}
For the inequality, apply \eqref{eq:spectral-block-bound} under Assumptions~\ref{ass:trace-admissibility}\ref{ass:diagonal-density}--\ref{ass:uniform-spectral-scale} with
$d=d_i$ and use $(m/n)^{2d_i}\le1$; the case $m=0$ is immediate.
For fixed $n$ and $m$, stationarity, boundedness of $n^{-D}$ and the
triangle inequality imply that the corresponding partial sums
converge in $L^2(\Omega;\mathbb H_2)$ as $L\to\infty$, since
$G_{\mathcal Q}^{(L)}\to G_{\mathcal Q}$ in
$L^2(\mathbb H_1,\gamma_Q;\mathbb H_2)$.
Passing to the limit in the displayed estimate and taking the
supremum over $m=\lfloor ns\rfloor$, $s\in[0,1]$, proves
\eqref{eq:projection-remainder-bound}.
\end{proof}

\section{Applications}
\label{sec:examples}

\subsection{Nonlinear transformations of a functional FARIMA process}
\label{sec:farima-application}

Fractionally integrated autoregressive moving average (FARIMA) models form a popular class exhibiting long-range dependence. The scalar construction goes back to \cite{GrangerJoyeux1980,Hosking1981}.  Functional FARIMA models
were introduced in \cite{LiRobinsonShang2020}. In this section, we apply Corollary~\ref{thm:common-hermite-limit} to
nonlinear transformations of a functional FARIMA process. These nonlinear transformations may include  
sample moment operators, sample eigenvalues or neural operators; see, for instance Section 7 of \cite{DukerZoub26}.

Let \(\Phi\) and \(\Theta\) be polynomials with real coefficients,
\(\Phi(0)=\Theta(0)=1\), and no zeros in \(\{z:|z|\le1\}\).
Fix \(d\in(0,1/2)\) and set \(\alpha=2d\). With respect to the basis
\(\{\phi_k\}_{k\ge1}\) of \(\mathbb H_1\) used above, let
\(\{\varepsilon_t\}_{t\in\mathbb Z}\) be independent centered Gaussian
random variables with covariance operator
\(\Sigma_\varepsilon\phi_k=\nu_k\phi_k\), where
\begin{equation}
 \nu_k>0\quad\text{for every }k\ge1,
 \qquad \sum_{k\ge1}\nu_k<\infty.
 \label{eq:farima-innovation-trace}
\end{equation}
We consider the stationary causal solution of
\begin{equation}
 \Phi(B)(1-B)^d X_t=\Theta(B)\varepsilon_t,
 \qquad t\in\mathbb Z,
 \label{eq:farima-functional-equation}
\end{equation}
where \(B\) denotes the backward-shift operator and the fractional
difference is defined by its binomial series. 

To identify the spectral quantities in
Assumption~\ref{ass:trace-admissibility}, define
\begin{equation}
 \begin{gathered}
 g_d(\omega)
 =\frac{1}{2\pi}
   \left|\frac{\Theta(e^{-i\omega})}{\Phi(e^{-i\omega})}\right|^2
   |1-e^{-i\omega}|^{-\alpha},
 \quad 
 0<|\omega|\le\pi,\\
 v_d=\int_{-\pi}^{\pi}g_d(\omega)\,\dd\omega,
 \quad
 \kappa=\frac{1}{2\pi v_d}
             \left|\frac{\Theta(1)}{\Phi(1)}\right|^2.
 \end{gathered}
 \label{eq:farima-common-spectral-shape}
\end{equation}
We verify below that \(0<v_d<\infty\). The covariance eigenvalues
are \(\lambda_k=v_d\nu_k\), and the standardized coordinate processes
\(\{W_{\phi_k}(X_t)\}_{t\in\mathbb Z}\), defined in
\eqref{eq:white-noise-mapping}, are independent copies of a scalar
Gaussian FARIMA process with variance one.
For a function $G\in L^2(\mathbb H_1,\gamma_{Q};\mathbb H_2)$, recall its Hermite coefficients, rank \(q\), and chaos
projections \(G_p\) from 
\eqref{eq:full-hermite-expansion}--\eqref{eq:chaotic-components}.

\begin{proposition}
\label{prop:farima-nonlinear-limit}
Let \(\{X_t\}_{ t \in \ZZ}\) satisfy
\eqref{eq:farima-functional-equation}, and let $G\in L^2(\mathbb H_1,\gamma_{Q};\mathbb H_2)$ be centered
and nonzero.
Suppose that \(q(1-\alpha)<1\), and let \(H\) be as in
\eqref{eq:common-H}. Then the interpolated partial-sum process
\(\widetilde S_n\) in \eqref{eq:continuous-interpolation}, with
\(T_n=n^{-H}I_{\mathbb H_2}\), satisfies
\begin{equation}
 \widetilde S_n\xrightarrow{d}Z_H^{(q)}
 \quad\text{in }C([0,1];\mathbb H_2),
 \qquad n\to\infty,
 \label{eq:farima-functional-limit}
\end{equation}
where \(Z_H^{(q)}\) is given in \eqref{eq:common-limit-process}
with \(\kappa_k=\kappa\) for every \(k\ge1\).
The process \(S_n\) in \eqref{eq:continuous-interpolation} converges
to the same limit in finite-dimensional distributions in \(\mathbb H_2\).
\end{proposition}

\begin{proof}
We obtain the convergence by means of
Corollary~\ref{thm:common-hermite-limit}. We need to verify
Assumptions \ref{ass:trace-admissibility} and
\ref{ass:common-active-exponent}.
The root conditions on $\Phi$ and $\Theta$ imply that
$\left|\Theta(e^{-i\omega})/\Phi(e^{-i\omega})\right|^2$
is bounded above and away from zero on $[-\pi,\pi]$. Moreover,
\[
 \frac{2}{\pi}|\omega|
 \le |1-e^{-i\omega}|
 \le |\omega|,
 \qquad |\omega|\le\pi.
\]
Consequently, $v_d$ in
\eqref{eq:farima-common-spectral-shape} satisfies
$0<v_d<\infty$, since $2d=\alpha\in(0,1)$.

Since $\Phi$ has no zeros in the closed unit disc, the rational
function $\Theta/\Phi$ is analytic on a neighborhood of that disc
and bounded there. The binomial series defining $(1-z)^{-d}$
converges for $|z|<1$. Hence the transfer function
\[
 A(z)\doteq\frac{\Theta(z)}{\Phi(z)}(1-z)^{-d}
      =\sum_{j=0}^{\infty}a_jz^j,
 \qquad |z|<1,
\]
is analytic on the open unit disc, with real coefficients $a_j$.

Since $0<d<1/2$ and $\Phi$ has no zeros in the closed unit
disc, the scalar FARIMA representation implies that $(a_j)_{j\ge0}$
is square summable and that
\[
 \sum_{j=0}^{\infty}a_j e^{-ij\omega}
 =
 \frac{\Theta(e^{-i\omega})}{\Phi(e^{-i\omega})}
 (1-e^{-i\omega})^{-d}
 \quad\text{in }L^2(-\pi,\pi);
\]
see Section~2.4.2 in \cite{PipirasTaqqu17}.
Consequently, Parseval's identity gives
\[
 \sum_{j=0}^{\infty}a_j^2
 =\frac{1}{2\pi}\int_{-\pi}^{\pi}
   \left|\frac{\Theta(e^{-i\omega})}
                {\Phi(e^{-i\omega})}\right|^2
   |1-e^{-i\omega}|^{-2d}\,\dd\omega
 =v_d<\infty.
\]

Since the innovations are independent and centered, the series
$X_t=\sum_{j\ge0}a_j\varepsilon_{t-j}$ converges in
$L^2(\Omega;\mathbb H_1)$, because
\[
 \sum_{j\ge0}
 \E\|a_j\varepsilon_{t-j}\|_{\mathbb H_1}^2
 =v_d\sum_{k\ge1}\nu_k<\infty.
\]
It defines a centered stationary Gaussian solution of
\eqref{eq:farima-functional-equation}, with
\begin{equation}
 \mathscr F_\omega=g_d(\omega)\Sigma_\varepsilon,
 \qquad Q=v_d\Sigma_\varepsilon,
 \qquad f_k(\omega)=\nu_k g_d(\omega).
 \label{eq:farima-diagonal-spectrum}
\end{equation}
In particular, $\lambda_k=v_d\nu_k$ for every $k\ge1$.
Since $\Sigma_\varepsilon$ is trace-class,
$\mathscr F_\omega$ is trace-class for almost every
$\omega\in[-\pi,\pi]$, and
\[
 \int_{-\pi}^{\pi}\Tr(\mathscr F_\omega)\,\dd\omega
 =v_d\sum_{k\ge1}\nu_k<\infty.
\]
This verifies
Assumption~\ref{ass:trace-admissibility}\ref{ass:trace-mathF}.
Furthermore, the spectral density and autocovariance operators
are diagonal in the basis $\{\phi_k\}_{k\ge1}$. Their spectral
densities $f_k$ are nonnegative and even, where evenness follows
from the real coefficients of $\Phi$ and $\Theta$.
Thus
Assumption~\ref{ass:trace-admissibility}\ref{ass:diagonal-density}
also holds.

For
Assumption~\ref{ass:trace-admissibility}\ref{ass:coordinatewise-lrd},
take $\alpha_k=\alpha$, $\kappa_k=\kappa$, and
\begin{equation}
 L_k(r)
 =\left|
   \frac{\Theta(e^{-ir})\Phi(1)}
        {\Phi(e^{-ir})\Theta(1)}
  \right|^2
  \left(\frac{r}{2\sin(r/2)}\right)^\alpha,
 \qquad 0<r\le\pi,\quad k\ge1.
 \label{eq:farima-slow-factor}
\end{equation}
Then \eqref{eq:global-spectral-envelope} follows from
\eqref{eq:farima-common-spectral-shape} and
\eqref{eq:farima-diagonal-spectrum}. The functions $L_k$ are
nonnegative, identical for all $k$, and extend continuously
to $[0,\pi]$, with value one at zero. Consequently,
$L_k(r)\to1$ as $r\downarrow0$, and there exists a finite
constant $M$ such that $L_k(r)\le M$ for every $k\ge1$
and $0<r\le\pi$. This verifies
Assumption~\ref{ass:trace-admissibility}\ref{ass:coordinatewise-lrd}.

To verify
Assumption~\ref{ass:trace-admissibility}\ref{ass:uniform-spectral-scale},
note that
\[
 \frac{f_k(\omega)}{\lambda_k}
 =\frac{\nu_k g_d(\omega)}{\nu_k v_d}
 =\frac{g_d(\omega)}{v_d},
 \qquad k\ge1.
\]
Thus all coordinates have the same standardized spectral density
and, in particular, the same spectral scale $\kappa_k=\kappa$.
Consequently,
\[
 \kappa_*=\kappa
 =\frac{1}{2\pi v_d}
   \left|\frac{\Theta(1)}{\Phi(1)}\right|^2<\infty,
\]
where finiteness follows from $v_d>0$ and $\Phi(1)\neq0$.
This completes the verification of
Assumption~\ref{ass:trace-admissibility}.

We turn to Assumption~\ref{ass:common-active-exponent}.
Square integrability of $G$ and its Hermite rank $q$ are
part of the hypotheses. Hence
Assumption~\ref{ass:common-active-exponent}\ref{ass:GinL2}
is satisfied.

For
Assumption~\ref{ass:common-active-exponent}\ref{ass:common-active-dom-chaos},
recall that $\alpha_j=\alpha$ for every input coordinate
$j\ge1$. In particular, whenever
$c_{i,\boldsymbol\ell}\neq0$ with
$\boldsymbol\ell\in\mathcal L_q$, we have
$\alpha_j=\alpha$ for every
$j\in\operatorname{supp}(\boldsymbol\ell)$. Moreover, the assumed
condition $q(1-\alpha)<1$ is equivalent to
$\alpha>1-1/q$. Thus
Assumption~\ref{ass:common-active-exponent}\ref{ass:common-active-dom-chaos}
holds.

Finally, by \eqref{eq:unequal-effective-exponent}, for every
$\boldsymbol\ell\in\mathcal L_p$ with $p>q$,
\[
 H(\boldsymbol\ell)
 =1-\frac p2(1-\alpha)
 <1-\frac q2(1-\alpha)
 =H,
\]
where the strict inequality follows from $\alpha<1$.
This verifies
Assumption~\ref{ass:common-active-exponent}\ref{ass:common-active-high-chaos}.
The conclusion follows from
Corollary~\ref{thm:common-hermite-limit}.
\end{proof}

For Hilbert-valued linear processes,
\cite{CharaciejusRackauskas2014,Duker2018} establish functional
central limit theorems with Gaussian limits and operator
normalizations. Their results concern partial sums of the linear
process itself and allow innovations that need not be Gaussian.
In the present FARIMA setting, taking \(G[x]=x\) gives
\(q=1\) and \(H=1/2+d\), consistent with the scalar normalization
in this linear limit theory. For general nonlinear transformations,
the normalization depends additionally on the Hermite rank, as
specified in \eqref{eq:common-H}.

A related functional FARIMA construction is given in
\cite[Section~4]{LiRobinsonShang2020}, where scalar fractional
integration is applied to a functional ARMA process. Their
Theorem~1 and Corollary~1 establish Gaussian limits for temporal
sums of curve time series through a decomposition into
finite-dimensional long-memory subspaces. In particular, their
FARIMA application imposes a finite-rank condition on the
covariance governing the dominant component. In our application,
infinitely many input coordinates may share the same memory
exponent, and we consider the partial-sum process of a nonlinear
transformation with values in \(\mathbb H_2\).

The input process considered here also belongs to the
operator-valued FIARMA framework of \cite{DurandRoueff2024}.
Indeed, we obtain \eqref{eq:farima-functional-equation} by taking
their fractional integration operator to be
\(dI_{\mathbb H_1}\), their autoregressive and moving-average
polynomials to be \(\Phi(z)I_{\mathbb H_1}\) and
\(\Theta(z)I_{\mathbb H_1}\), respectively, and their innovations
to be Gaussian. Their framework permits more general operator
coefficients and provides existence results for normal memory
operators, together with prediction and estimation theory.
Our proposition establishes a nonlinear functional limit theorem
for the Gaussian subclass specified above.

Hilbert-valued nonlinear non-central limit theorems have also
been obtained for scalar Gaussian inputs. In
\cite[Theorem~2.7 of version~3]{BetkenDehling2021}, a scalar
long-memory Gaussian sequence is transformed into an observation
with values in a weighted \(L^2\) space. Their reduction principle
identifies the leading Hermite term and yields convergence of
the normalized sum at time one. The limit is a deterministic
coefficient function multiplied by a scalar Hermite random
variable and therefore lies in a single deterministic output
direction. Thus infinite-dimensional output is already covered
by that result. In the present setting, the input may also be
infinite-dimensional, and mixed Hermite terms involving different
input coordinates may survive in the limit. Consequently, the
limiting covariance need not have finite rank. Moreover, under
our FARIMA assumptions, we obtain convergence of the interpolated
partial-sum process in \(C([0,1];\mathbb H_2)\).

\subsection{Excursion areas of spherical Gaussian fields}
\label{subsec:spherical-excursions}

We return to the sojourn functionals that motivated this work and apply
Theorem~\ref{thm:operator-normalized-full} to a spherical Gaussian field
sampled at integer times. We first identify the dominant Hermite indices
of the excursion area, then state the resulting non-central limit theorem,
and finally compare it with the continuous-time results of
Marinucci, Rossi and Vidotto~\cite{MarinucciRV}.

\subsubsection{The spherical model and the excursion functional}
\label{subsec:spherical-setting}

We use the spherical random wave setting presented in 
\cite{MarinucciRV}. Let
$\{Y_{\ell,m}:\ell\ge0,\ |m|\le\ell\}$ be the standard real orthonormal basis of
$\mbbH_1=L^2(\mathbb S^2,dx)$ consisting of spherical harmonics, where
$dx$ is surface measure with total mass $4\pi$ and
$Y_{0,0}=(4\pi)^{-1/2}$. Let
$\{Z(x,t):x\in\mathbb S^2,\ t\in\RR\}$ be a jointly measurable,
centered Gaussian field, stationary in time and isotropic in space,
admitting the representation
\begin{equation}
Z(x,t)=\sum_{\ell\ge0}\sum_{m=-\ell}^{\ell}
a_{\ell,m}(t)Y_{\ell,m}(x), \qquad 
\E[a_{\ell,m}(t)a_{\ell',m'}(s)] =\mathds{1}_{\{\ell=\ell',\,m=m'\}}C_\ell(t-s),
\label{eq:spherical-field}
\end{equation}
with convergence in $L^2(\Omega\times\mathbb S^2)$ for every $t$.
Each multipole
$Z_\ell(\cdot,t)=\sum_m a_{\ell,m}(t)Y_{\ell,m}$ is a
random eigenfunction of the spherical Laplacian with eigenvalue
$-\ell(\ell+1)$, and we assume that $C_\ell(0)>0$ and
\begin{equation}
 \sum_{\ell\ge0}\frac{2\ell+1}{4\pi}C_\ell(0)=1,
 \label{eq:spherical-unit-variance}
\end{equation}
so that $\E Z(x,t)^2=1$ by \eqref{eq:spherical-field}. Hence the sampled sequence
$X_t \doteq Z(\cdot,t)$, $t\in\ZZ$, is centered, stationary and Gaussian
with values in $\mbbH_1$, and its covariance operator satisfies
$QY_{\ell,m}=C_\ell(0)Y_{\ell,m}$ and $\Tr(Q)=4\pi$.

Throughout this section, suppose that $\{X_t\}_{t\in\ZZ}$ satisfies
Assumption~\ref{ass:trace-admissibility} in the spherical harmonic basis.
Index the input coordinates by
$\mathcal J=\{(\ell,m):\ell\ge0,\ |m|\le\ell\}$ through a fixed
enumeration, and write $\ell(j)=\ell$, $Y_j=Y_{\ell,m}$ and
$\lambda_j=C_\ell(0)$ for $j=(\ell,m)$. The standardized coordinate
processes $\xi_j(t)=a_{\ell,m}(t)/\sqrt{C_\ell(0)}$ are independent.
Since their laws depend only on $\ell$, we write the parameters in
Assumption~\ref{ass:trace-admissibility} as
$\alpha_j=1-\beta_{\ell(j)}$ and $\kappa_j=\kappa_{\ell(j)}$.
Thus $\beta_\ell\in(0,1)$ is the memory parameter of the degree-$\ell$
coefficient processes, and $\kappa_\ell>0$ is their standardized
spectral scale, so that, for $j=(\ell,m)$,
\[
 \frac{f_j(\omega)}{C_\ell(0)}
 \sim\kappa_\ell|\omega|^{\beta_\ell-1},
 \qquad \omega\to0.
\]
In particular, $\beta_0=1-\alpha_{(0,0)}$ is the memory parameter of
$a_{0,0}$, the coefficient of the constant spherical harmonic, determining the monopole $Z_0(x,t)=a_{0,0}(t)/\sqrt{4\pi}$.
Smaller values of $\beta_\ell$ correspond to stronger memory. Define
\begin{equation}
 \beta_*\doteq\inf_{\ell\ge1}\beta_\ell,\qquad
 \mathcal I^*\doteq\{\ell\ge1:\beta_\ell=\beta_*\}.
 \label{eq:spherical-strongest-memory}
\end{equation}
From Assumption~\ref{ass:trace-admissibility}, we have that
$\inf_{\ell\ge0}\beta_\ell=\min\{\beta_0,\beta_*\}>0$.
Denote by $\mathbf e_j$ the unit multiplicity index at $j\in\mathcal J$.

Let \(\varphi\) and \(\Phi\) denote the standard Gaussian density
and distribution function, respectively. For \(u\in\RR\), define
\begin{equation}
\label{eq:spherical-functional}
 G_u[f]
 =
 \int_{\mathbb S^2}\mathds{1}_{\{f(x)\ge u\}}\,dx,
 \qquad
 \overline G_u[f]
 =
 G_u[f]-4\pi(1-\Phi(u)),
 \qquad f\in\mbbH_1.
\end{equation}
For a Hermite multi-index $\boldsymbol r\in\mathcal L_p$, let
$j_a=(\ell_a,m_a)$, $1\le a\le p$, list its input coordinates, with
each $j\in\mathcal J$ appearing exactly $r_j$ times. Define the
\textit{generalized Gaunt integrals}
\begin{equation} \label{eq:Gaunt-integrals}
     \mathcal G_{\ell_1\cdots\ell_p}^{m_1\cdots m_p}
 =
 \int_{\mathbb S^2}
 \prod_{a=1}^pY_{\ell_a,m_a}(x)\,dx.
\end{equation}
The following lemma identifies the Hermite expansion of
\(G=\overline G_u:\mbbH_1\to\RR\).
\begin{lemma} \label{lem:spherical-hermite-decomposition}
For every \(u\in\RR\) and \(t\in\ZZ\),
\begin{equation}
\label{eq:spherical-hermite-coefficients}
\overline G_u[X_t]=\sum_{p=1}^{\infty}\sum_{\boldsymbol r\in\mathcal L_p}
 c_{u,\boldsymbol r}H_{\boldsymbol r}[X_t], \qquad c_{u,\boldsymbol r}
 =\frac{\varphi(u)H_{p-1}(u)}{\boldsymbol r!}
 \left(\prod_{a=1}^p C_{\ell_a}(0)^{1/2}\right)
 \mathcal G_{\ell_1\cdots\ell_p}^{m_1\cdots m_p},
\end{equation}
for $\boldsymbol r\in\mathcal L_p$, where the series converges in $L^2(\Omega)$ and $H_{\boldsymbol r}[X_t]
=\prod_{j\in\mathcal J}H_{r_j}(\xi_j(t))$.
\end{lemma}

\begin{proof}
From Lemma~4.1 of \cite{MarinucciRV}, we can write
\begin{align}
 \overline G_u[Z(\cdot,t)]
 &=
 \sum_{p=1}^{\infty}
 \frac{\varphi(u)H_{p-1}(u)}{p!}
 \int_{\mathbb S^2}H_p(Z(x,t))\,dx\\
 &=
 \sum_{p=1}^{\infty}
 \frac{\varphi(u)H_{p-1}(u)}{p!}
 \int_{\mathbb S^2}
 \sum_{\boldsymbol r\in\mathcal L_p}
 \frac{p!}{\boldsymbol r!}
 \left(
 \prod_{a=1}^p
 C_{\ell_a}(0)^{1/2}Y_{\ell_a,m_a}(x)
 \right)
 H_{\boldsymbol r}[X_t]\,dx,
\end{align}
where the second line follows from decomposing
\(Z(x,t)=\sum_j C_{\ell(j)}(0)^{1/2}Y_j(x)\xi_j(t)\), the fact that
\(\sum_j C_{\ell(j)}(0)Y_j(x)^2=1\) by
\eqref{eq:spherical-unit-variance}, and the multinomial Hermite
formula. 
\end{proof}

Recall from \eqref{eq:unequal-effective-exponent} that
\[
 \beta(\boldsymbol r)=\sum_{a=1}^p\beta_{\ell_a},
 \qquad H(\boldsymbol r)=1-\tfrac12\beta(\boldsymbol r).
\]
Thus the active terms with the smallest $\beta(\boldsymbol r)$ have
the largest effective exponent $H(\boldsymbol r)$ and determine the
normalization. To describe the possible dominant cubic terms, define
\[
 \mathcal C_*
 =
 \left\{
 \boldsymbol r\in\mathcal L_3:
 \ell_1,\ell_2,\ell_3\in \mathcal I^*,
 \quad
 \mathcal G_{\ell_1\ell_2\ell_3}^{m_1m_2m_3}\ne0
 \right\}.
\]

We call $u$ \emph{admissible} if either (i) $u\ne0$ and
$\beta_0\le2\beta_*$, (ii) $u\ne0$ and $\mathcal I^*\ne\varnothing$,
(iii) $u=0$ and $\beta_0\le3\beta_*$, or (iv) $u=0$ and
$\mathcal I^*$ contains an even integer. Under these four cases, the smallest
effective memory parameter is attained by an active Hermite term.
The following lemma identifies the corresponding exponents and
dominant sets.

\begin{lemma}
\label{lem:spherical-dominant}
In the setting above, let \(u\) be admissible. Then
\(\overline G_u\) is centered, bounded and nonzero, and satisfies
Assumption~\ref{ass:bigAssump2}.
For \(u\ne0\), its exponent and dominant set are
\begin{equation}
\label{eq:spherical-dominant-nonzero}
 d_u=1-\tfrac12\min\{\beta_0,2\beta_*\},\quad 
 \mathcal D_u=
 \begin{cases}
 \{\mathbf e_{(0,0)}\},
     &\beta_0<2\beta_*,\\
 \{2\mathbf e_{(\ell,m)}:\ell\in \mathcal I^*,\ |m|\le\ell\},
     &2\beta_*<\beta_0,\\
 \{\mathbf e_{(0,0)}\}
 \cup\{2\mathbf e_{(\ell,m)}:\ell\in \mathcal I^*,\ |m|\le\ell\},
     &\beta_0=2\beta_*.
 \end{cases}
\end{equation}
For \(u=0\), they are
\begin{equation}
\label{eq:spherical-dominant-zero}
 d_0=1-\tfrac12\min\{\beta_0,3\beta_*\},\quad 
 \mathcal D_0=
 \begin{cases}
 \{\mathbf e_{(0,0)}\},&\beta_0<3\beta_*,\\
 \mathcal C_*,&3\beta_*<\beta_0,\\
 \{\mathbf e_{(0,0)}\}\cup\mathcal C_*,
     &\beta_0=3\beta_*.
 \end{cases}
\end{equation}
Moreover, \(\mathcal C_*\ne\varnothing\) if and only if
\(\mathcal I^*\) contains an even integer.
\end{lemma}

\begin{proof}
First, $|\overline G_u|\le4\pi$ and $\E\overline G_u[X_0]=0$. We identify the dominant terms by comparing their effective exponents
$H(\boldsymbol r)=1-\beta(\boldsymbol r)/2$ in the different parameter
regimes. We then verify Assumption~\ref{ass:bigAssump2} and explain why
the remaining chaoses are negligible under the resulting normalization.
By orthonormality,
\begin{equation}
\label{eq:Gaunt-first-second}
 \mathcal G_\ell^m
 =\sqrt{4\pi}\,\mathds{1}_{\{(\ell,m)=(0,0)\}},
 \qquad
 \mathcal G_{\ell_1\ell_2}^{m_1m_2}
 =\mathds{1}_{\{\ell_1=\ell_2,\ m_1=m_2\}},
\end{equation}
implying that the only active first-order index is
$\mathbf e_{(0,0)}$, and the only possible second-order indices are
$2\mathbf e_{(\ell,m)}$, with coefficients
\begin{equation}
\label{eq:spherical-low-coefficients}
 c_{u,\mathbf e_{(0,0)}}
 =\varphi(u)\sqrt{4\pi C_0(0)},
 \qquad
 c_{u,2\mathbf e_{(\ell,m)}}
 =\frac{u\varphi(u)C_\ell(0)}2.
\end{equation}
The first coefficient is nonzero for every $u$, so $\overline G_u$ is
nonzero and has an active term with effective exponent $1-\beta_0/2$.
Moreover, any index $\boldsymbol r$ of order $p\ge2$ involving the
monopole satisfies $\beta(\boldsymbol r)
 \ge\beta_0+(p-1)\min\{\beta_0,\beta_*\}>\beta_0$. Its effective exponent is therefore strictly smaller than that of the
linear monopole term, in every parameter regime.

We first consider $u\ne0$. All the quadratic coefficients in
\eqref{eq:spherical-low-coefficients} are nonzero. For $\ell\ge1$, we have that $H(2\mathbf e_{(\ell,m)})=1-\beta_\ell\le1-\beta_*$, and equality holds when $\ell\in\mathcal I^*$. Thus the
competition is between the monopole exponent $1-\beta_0/2$ and the
supremum $1-\beta_*$ of the quadratic exponents over positive degrees.
If $\beta_0<2\beta_*$, the monopole exponent is larger. If
$2\beta_*<\beta_0$, the quadratic exponent is larger and attained on $\mathcal I^*\ne\varnothing$.
If $\beta_0=2\beta_*$, the monopole and all quadratic terms with
$\ell\in\mathcal I^*$ have the same effective exponent and must be
retained together; if $\mathcal I^*=\varnothing$, only the monopole
attains it.

To exclude higher chaoses, it remains to consider indices that do not
involve the monopole. For every such index of order $p\ge3$, $ \beta(\boldsymbol r)\ge p\beta_*\ge3\beta_*>2\beta_*$, meaning $H(\boldsymbol r)<1-\beta_*\le d_u$. Together with the preceding bound for indices involving the monopole,
this shows that no higher-order term can attain the dominant exponent.
The three comparisons above therefore give
\eqref{eq:spherical-dominant-nonzero}.

We consider $u=0$. Every even chaos vanishes because
$H_{p-1}(0)=0$ for even $p$, whereas the first-order coefficient remains
nonzero. Since $H_2(0)=-1$, a cubic coefficient is nonzero exactly when
$G_{\ell_1\ell_2\ell_3}^{m_1m_2m_3}\ne0$. For cubic indices not involving the
monopole, $ H(\boldsymbol r)
 =1-\tfrac12(\beta_{\ell_1}+\beta_{\ell_2}+\beta_{\ell_3})
 \le1-\tfrac32\beta_*$, and equality holds precisely for the active indices in $\mathcal C_*$.
To determine when this set is nonempty, recall that
$Y_{\ell,m}(-x)=(-1)^\ell Y_{\ell,m}(x)$. Hence a Gaunt integral
vanishes whenever the sum of its degrees is odd. In particular,
$\mathcal C_*=\varnothing$ if every degree in $\mathcal I^*$ is odd.
The converse is also true, since for every even degree $\ell$, the standard zonal harmonic
satisfies $\mathcal G_{\ell\ell\ell}^{000}>0$, so
$\mathcal C_*\ne\varnothing$ whenever $\mathcal I^*$ contains an even
integer.

We therefore compare the monopole exponent $1-\beta_0/2$ with the
cubic bound $1-3\beta_*/2$. If $\beta_0<3\beta_*$, the monopole
dominates every cubic term. If $3\beta_*<\beta_0$, admissibility
ensures that $\mathcal C_*\ne\varnothing$, and precisely these cubic
terms attain the larger exponent. If $\beta_0=3\beta_*$, the monopole
and the terms in $\mathcal C_*$ have equal effective exponents, with
only the monopole retained when $\mathcal C_*=\varnothing$.
Every remaining active cubic index has a smaller effective exponent
than $d_0$. Finally, for indices of odd order $p\ge5$ not involving
the monopole, $\beta(\boldsymbol r)\ge p\beta_*\ge5\beta_*>3\beta_*$, so that $H(\boldsymbol r)<1-\tfrac32\beta_*\le d_0$ proving \eqref{eq:spherical-dominant-zero}.

In both cases, the smallest effective memory parameter is positive
and at most $\beta_0<1$, so $1/2<d_u<1$. We can then conclude that
$\mathcal D_u\ne\varnothing$, and, moreover,
\[
 \frac{\E|\overline G_u[X_0]|^2}{2d_u-1}
 \le\frac{16\pi^2}{2d_u-1}<\infty,
\]
so Assumption~\ref{ass:bigAssump2} is satisfied. 
\end{proof}

\subsubsection{The non-central limit theorem}
\label{subsec:spherical-limit}

Before stating the result, we describe the limiting processes, whose
values at $s=1$ give the limiting random variables for the normalized
excursion sums. On the common spectral probability space of
Section~\ref{subse-Wiender}, for $\boldsymbol r\ne0$ with
$\beta(\boldsymbol r)<1$ and $s \in [0,1]$, set
\begin{equation}
 \mathcal Z_{\boldsymbol r}(s)
 =\kappa^{\boldsymbol r/2}
 I_{|\boldsymbol r|}(J_{\boldsymbol r}K_{s,\boldsymbol r}),
 \qquad
 \kappa^{\boldsymbol r/2}=\prod_j\kappa_{\ell(j)}^{r_j/2},
 \label{eq:spherical-building-block}
\end{equation}
where $J_{\boldsymbol r}$ and $K_{s,\boldsymbol r}$ are defined in
\eqref{eq:labelled-lift} and \eqref{eq:operator-limit-kernel},
respectively. Define
\begin{align}
 \mathcal B_u(s)
 &=\varphi(u)\sqrt{4\pi C_0(0)}\,
 \mathcal Z_{\mathbf e_{(0,0)}}(s),
 \label{eq:spherical-B}\\
 \mathcal R_u(s)
 &=\frac{u\varphi(u)}2\sum_{\ell\in \mathcal I^*}C_\ell(0)
 \sum_{m=-\ell}^{\ell}\mathcal Z_{2\mathbf e_{(\ell,m)}}(s),
 \qquad 2\beta_*<1,
 \label{eq:spherical-R}\\
 \mathcal T(s)
 &=-\varphi(0)\sum_{\boldsymbol r\in\mathcal C_*}
 \frac{\mathcal G_{\ell_1\ell_2\ell_3}^{m_1m_2m_3}}{\boldsymbol r!}
 \left(\prod_{a=1}^3 C_{\ell_a}(0)^{1/2}\right)
 \mathcal Z_{\boldsymbol r}(s),
 \qquad 3\beta_*<1,
 \label{eq:spherical-T}
\end{align}
with the convention that the sum over an empty set vanishes. These series converge in
$L^2(\Omega)$ for every $s$, by the labeled isometry
\eqref{eq:labelled-isometry}, the uniform kernel bound
\eqref{eq:uniform-kernel-norm}, and
$\sum_{\boldsymbol r}\boldsymbol r!|c_{u,\boldsymbol r}|^2
=\E|\overline G_u[X_0]|^2<\infty$.

The process $\mathcal B_u$ is a multiple of a fractional Brownian motion
with Hurst index $1-\beta_0/2$; the process $\mathcal R_u$ is a weighted sum of
independent Rosenblatt processes with Hurst index $1-\beta_*$, which
follows from their first- and second-order spectral kernels, respectively.
If $\mathcal C_*\ne\varnothing$, then $\mathcal T$ is a nonzero process
in the third Wiener chaos.
Moreover, $\mathcal B_u$ depends only on $\widehat W_{(0,0)}$, whereas
$\mathcal R_u$ and $\mathcal T$ depend only on
$\{\widehat W_{(\ell,m)}:\ell\in\mathcal I^*,\ |m|\le\ell\}$, so that $\mathcal B_u$ is independent
of the latter processes.

For an admissible threshold, define
\begin{equation}
 \mathcal E_u=
 \begin{cases}
 \mathcal B_u, & \beta_0<2\beta_*,\\
 \mathcal R_u, & 2\beta_*<\beta_0,\\
 \mathcal B_u+\mathcal R_u, & \beta_0=2\beta_*,
 \end{cases}
\;  u\ne0, \qquad 
 \mathcal E_0=
 \begin{cases}
 \mathcal B_0, & \beta_0<3\beta_*,\\
 \mathcal T, & 3\beta_*<\beta_0,\\
 \mathcal B_0+\mathcal T, & \beta_0=3\beta_*.
 \end{cases}
 \label{eq:spherical-limits}
\end{equation}
The conditions for defining $\mathcal R_u$ and $\mathcal T$ hold whenever
these processes occur in \eqref{eq:spherical-limits}, since $\beta_0<1$. In the equality cases, the two components are independent;
the limit is non-Gaussian whenever the higher-order component is nonzero.

\begin{theorem}
\label{prop:spherical-excursions}
Let $X_t=Z(\cdot,t)$ be the spherical Gaussian sequence above, satisfying
Assumption~\ref{ass:trace-admissibility} in the spherical harmonic basis.
Let $u_1,\ldots,u_r$ be admissible thresholds, and set
$G=(\overline G_{u_1},\ldots,\overline G_{u_r})$ and
$D=\operatorname{diag}(d_{u_1},\ldots,d_{u_r})$, with the exponents of
Lemma~\ref{lem:spherical-dominant}. Then, as $n\to\infty$,
\[
 n^{-D}\sum_{t=1}^{\lfloor ns\rfloor}G[X_t]
 \xrightarrow{\operatorname{f.d.d.}}
 (\mathcal E_{u_1}(s),\ldots,\mathcal E_{u_r}(s)).
\]
\end{theorem}

\begin{proof}
By Lemma~\ref{lem:spherical-dominant}, each coordinate satisfies
item~\ref{ass:operator-standing} of Assumption~\ref{ass:bigAssump2}, and
\[
 \sum_{i=1}^r\frac{\E|\overline G_{u_i}[X_0]|^2}{2d_{u_i}-1}<\infty,
 \qquad \min_{1\le i\le r}d_{u_i}>\tfrac12.
\]
Theorem~\ref{thm:operator-normalized-full} therefore yields the
conclusion. Its $i$th limit coordinate is
$\sum_{\boldsymbol r\in\mathcal D_{u_i}}
c_{u_i,\boldsymbol r}\mathcal Z_{\boldsymbol r}$, which equals
$\mathcal E_{u_i}$ by \eqref{eq:spherical-hermite-coefficients} and the
dominant sets in Lemma~\ref{lem:spherical-dominant}.
\end{proof}

We present a couple of remarks regarding our theorem
\begin{remark}
\label{rem:spherical-comments}
If $u\ne0$, $2\beta_*<\beta_0$ and $\mathcal I^*=\varnothing$, the supremum
defining $d_u$ is not attained, so our main theorem does not apply.
For $u=0$ outside the admissible regimes, the active indices must be
determined from \eqref{eq:spherical-hermite-coefficients}. Whenever
$\min_{c_{0,\boldsymbol r}\ne0}\beta(\boldsymbol r)$ is attained,
the main theorem applies with $d_0$ equal to one minus half this minimum.
\end{remark}

The next remark compares our result with that of \cite{MarinucciRV}.

\begin{remark}
\label{rem:spherical-comparison}
We compare our result with those of \cite{MarinucciRV}, first
at the level of conclusions and then in terms of the assumptions
on the process \(\{X_t\}_{t \in \ZZ}\).

In the notation of \cite{MarinucciRV},
\(\mathcal A_u(t)=G_u[Z(\cdot,t)]\), and their asymptotic results
concern the continuous-time integral
\[
 \mathcal M_T(u)
 =
 \int_0^T
 \bigl(\mathcal A_u(t)-\E\mathcal A_u(t)\bigr)\,dt
 =
 \int_0^T\overline G_u[Z(\cdot,t)]\,dt.
\]
Their limits are stated after variance normalization, whereas our
normalization uses powers \(n^{-d_u}\). Their Theorems~2.2 and~2.7
establish convergence to Gaussian and composite Rosenblatt random
variables, respectively, with the distributional conclusion of
Theorem~2.7 assuming finitely many dominant multipoles.
Their Theorem~2.9 establishes a reduction to the third-chaos
component for \(u=0\), but does not establish its limiting
distribution; see Remark~2.10 therein.

Our Theorem~\ref{prop:spherical-excursions} concerns sampled
partial-sum processes and provides discrete-time counterparts
of these long-memory results. Passing from sampled sums to
continuous-time integrals would require a separate approximation
argument, which we do not provide here. Within the corresponding long-memory regimes, our
theorem provides the following additional conclusions:
(i) we identify the third-chaos process \(\mathcal T\) explicitly
and prove convergence to it;
(ii) we do not require finiteness of \(\mathcal I^*\) or
\(\mathcal C_*\);
(iii) we obtain joint finite-dimensional convergence at finitely
many thresholds \(u_1,\ldots,u_r\), allowing different normalization
exponents; and
(iv) we treat the equality cases \(\beta_0=2\beta_*\) and
\(\beta_0=3\beta_*\), discussed in Remark~2.11 of
\cite{MarinucciRV}.

The assumptions in \cite{MarinucciRV} are formulated in the time
domain and are therefore not immediately comparable to our
spectral Assumption~\ref{ass:trace-admissibility}. One difference is that
our assumption requires every retained multipole to have long
memory, so the short-memory and critical logarithmic regimes
considered in \cite{MarinucciRV} are not covered by our theorem.

Nevertheless, Remark~\ref{rem:time-domain-sufficient} allows us
to verify Assumption~\ref{ass:trace-admissibility} for a subclass
of their covariance family. Within the setting presented in Condition~1.3 of 
\cite{MarinucciRV}, suppose that
\begin{equation}
\label{eq:spherical-power-covariance}
 C_\ell(\tau)
 =
 G_\ell(\tau)(1+|\tau|)^{-\beta_\ell},
 \qquad
 0<\inf_{\ell\ge0}\beta_\ell
 \le\sup_{\ell\ge0}\beta_\ell<1,
\end{equation}
where \(C_\ell(0)>0\), and additionally assume that $ \sup_{\ell\ge0}\sum_{h=1}^{\infty}
 \left|\frac{G_\ell(h)}{C_\ell(0)}-1\right|
 <\infty$. These conditions are stronger than Condition~1.3 of 
\cite{MarinucciRV} since $ \sup_\ell \beta_\ell < 1$ and (ii) they require the absolute summability of $\left|\frac{G_\ell(h)}{C_\ell(0)}-1\right|$. Under these conditions, the conditions of Remark~\ref{rem:time-domain-sufficient} hold, and so
Assumption~\ref{ass:trace-admissibility} follows. We stress that these additional conditions are sufficient, rather than
necessary, for our spectral assumption.
\end{remark}

\small
\bibliographystyle{acm}
\bibliography{references}

@article{Breuer83,
  author = {Breuer, P. and Major, P.},
  title = {Central limit theorems for non-linear functionals of {Gaussian} fields},
  journal = {Journal of Multivariate Analysis},
  volume = {13},
  number = {3},
  pages = {425--441},
  year = {1983},
  doi = {10.1016/0047-259x(83)90019-2}
}

@article{Bulinski,
  author = {Bulinski, A. and Spodarev, E. and Timmermann, F.},
  title = {Central limit theorems for the excursion set volumes of weakly dependent random fields},
  journal = {Bernoulli},
  volume = {18},
  number = {1},
  pages = {100--118},
  year = {2012},
  doi = {10.3150/10-bej339}
}

@article{CampeseNourdin20,
  author = {Campese, S. and Nourdin, I. and Nualart, D.},
  title = {Continuous {Breuer--Major} theorem: Tightness and nonstationarity},
  journal = {The Annals of Probability},
  volume = {48},
  number = {1},
  pages = {147--177},
  year = {2020},
  doi = {10.1214/19-aop1357}
}

@article{Caponera,
  author = {Caponera, A. and Rossi, M. and Ruiz-Medina, M. D.},
  title = {Sojourn functionals of time-dependent {$\chi^2$}-random fields on two-point homogeneous spaces},
  journal = {Journal of Applied Probability},
  volume = {62},
  number = {4},
  pages = {1493--1512},
  year = {2025},
  doi = {10.1017/jpr.2025.18}
}

@article{Dobrushin,
  author = {Dobrushin, R. L. and Major, P.},
  title = {Non-central limit theorems for non-linear functionals of {Gaussian} fields},
  journal = {Zeitschrift f{\"u}r Wahrscheinlichkeitstheorie und Verwandte Gebiete},
  volume = {50},
  number = {1},
  pages = {27--52},
  year = {1979},
  doi = {10.1007/bf00535673}
}

@misc{DukerZoub25,
  author = {D{\"u}ker, M.-Ch. and Zoubouloglou, P.},
  title = {The fourth-moment theorem on {Hilbert} spaces},
  year = {2025},
  howpublished = {\href{https://arxiv.org/abs/2509.20342}{arXiv:2509.20342v3}},
  eprint = {2509.20342},
  archivePrefix = {arXiv},
  primaryClass = {math.PR},
  doi = {10.48550/arXiv.2509.20342},
  url = {https://arxiv.org/abs/2509.20342},
  note = {First posted in 2025; revised 22 September 2026 (v3)}
}

@article{DehlingTaqqu89,
  author = {Herold Dehling and Murad S. Taqqu},
  title = {The Empirical Process of some Long-Range Dependent Sequences with an Application to {U}-Statistics},
  journal = {The Annals of Statistics},
  volume = {17},
  number = {4},
  pages = {1767--1783},
  year = {1989},
  publisher = {Institute of Mathematical Statistics},
  doi = {10.1214/aos/1176347394}
}

@article{KraetschmerUrusov2023,
  author = {Kr{\"a}tschmer, Volker and Urusov, Mikhail},
  title = {A {Kolmogorov--Chentsov} Type Theorem on General Metric Spaces with Applications to Limit Theorems for {Banach}-Valued Processes},
  journal = {Journal of Theoretical Probability},
  volume = {36},
  number = {3},
  pages = {1454--1486},
  year = {2023},
  doi = {10.1007/s10959-022-01207-8}
}

@book{kallenberg2002,
  author = {Kallenberg, Olav},
  title = {Foundations of modern probability},
  edition = {Second},
  year = {2002},
  publisher = {Springer},
  doi = {10.1007/978-1-4757-4015-8}
}

@article{bai2019sensitivity,
  author = {Bai, Shuyang and Taqqu, Murad S},
  title = {Sensitivity of the {H}ermite rank},
  journal = {Stochastic Processes and their Applications},
  volume = {129},
  number = {3},
  pages = {822--840},
  year = {2019},
  publisher = {Elsevier},
  doi = {10.1016/j.spa.2018.03.020}
}

@book{billingsley1999convergence,
  author = {Billingsley, Patrick},
  title = {Convergence of probability measures},
  edition = {Second},
  year = {1999},
  publisher = {John Wiley \& Sons},
  address = {New York},
  doi = {10.1002/9780470316962}
}

@misc{DukerZoub26,
  author = {D{\"u}ker, M.-Ch. and Zoubouloglou, P.},
  title = {{Breuer--Major} theorems for {Hilbert} space-valued random variables},
  year = {2026},
  howpublished = {\href{https://arxiv.org/abs/2405.11452v3}{arXiv:2405.11452v3}},
  eprint = {2405.11452v3},
  archivePrefix = {arXiv},
  primaryClass = {math.PR},
  doi = {10.48550/arXiv.2405.11452},
  url = {https://arxiv.org/abs/2405.11452v3},
  note = {Originally posted in 2024; version 3 revised 22 September 2026}
}

@book{Ivanov89,
  author = {Ivanov, A. V. and Leonenko, N. N.},
  title = {Statistical Analysis of Random Fields},
  series = {Mathematics and its Applications (Soviet Series)},
  volume = {28},
  year = {1989},
  publisher = {Kluwer Academic Publishers},
  address = {Dordrecht, The Netherlands},
  doi = {10.1007/978-94-009-1183-3}
}

@book{Leonenko99,
  author = {Leonenko, N. N.},
  title = {Limit Theorems for Random Fields with Singular Spectrum},
  series = {Mathematics and its Applications},
  volume = {465},
  year = {1999},
  publisher = {Kluwer Academic Publishers},
  address = {Dordrecht},
  doi = {10.1007/978-94-011-4607-4}
}

@article{LeonenkoOlenko14,
  author = {Leonenko, N. N. and Olenko, A.},
  title = {Sojourn measures of {Student} and {Fisher--Snedecor} random fields},
  journal = {Bernoulli},
  volume = {20},
  number = {3},
  pages = {1454--1483},
  year = {2014},
  doi = {10.3150/13-bej529}
}

@article{GrangerJoyeux1980,
  author = {Granger, C. W. J. and Joyeux, Roselyne},
  title = {An introduction to long-memory time series models and fractional differencing},
  journal = {Journal of Time Series Analysis},
  volume = {1},
  number = {1},
  pages = {15--29},
  year = {1980},
  doi = {10.1111/j.1467-9892.1980.tb00297.x}
}

@article{Hosking1981,
  author = {Hosking, J. R. M.},
  title = {Fractional differencing},
  journal = {Biometrika},
  volume = {68},
  number = {1},
  pages = {165--176},
  year = {1981},
  doi = {10.1093/biomet/68.1.165}
}

@article{DurandRoueff2024,
  author = {Durand, Amaury and Roueff, Fran{\c{c}}ois},
  title = {{Hilbert} space-valued fractionally integrated autoregressive moving average processes with long memory operators},
  journal = {Journal of Statistical Planning and Inference},
  volume = {231},
  pages = {106146},
  year = {2024},
  doi = {10.1016/j.jspi.2024.106146}
}

@article{Leonenkoetal17,
  author = {Leonenko, N. N. and Ruiz-Medina, M. D. and Taqqu, M. S.},
  title = {Non-central limit theorems for random fields subordinated to gamma-correlated random fields},
  journal = {Bernoulli},
  volume = {23},
  number = {4B},
  pages = {3469--3507},
  year = {2017},
  doi = {10.3150/16-bej853}
}

@article{Leonenkoetal17b,
  author = {Leonenko, N. N. and Ruiz-Medina, M. D. and Taqqu, M. S.},
  title = {{Rosenblatt} distribution subordinated to {Gaussian} random fields with long-range dependence},
  journal = {Stochastic Analysis and Applications},
  volume = {35},
  number = {1},
  pages = {144--177},
  year = {2017},
  doi = {10.1080/07362994.2016.1230723}
}

@article{LeonenkoRuiz-Medina2023,
  author = {Leonenko, N. N. and Ruiz-Medina, M. D.},
  title = {Sojourn functionals for spatiotemporal {Gaussian} random fields with long memory},
  journal = {Journal of Applied Probability},
  volume = {60},
  number = {1},
  pages = {148--165},
  year = {2023},
  doi = {10.1017/jpr.2022.30}
}

@article{Leonenkoruizmedi25,
  author = {Leonenko, N. N. and Ruiz-Medina, M. D.},
  title = {High-level moving excursions for spatiotemporal {Gaussian} random fields with long range dependence},
  journal = {Journal of Statistical Physics},
  volume = {192},
  number = {2},
  pages = {19},
  year = {2025},
  doi = {10.1007/s10955-025-03396-y}
}

@article{Makogin20,
  author = {Makogin, V. and Spodarev, E.},
  title = {Limit theorems for excursion sets of subordinated {Gaussian} random fields with long-range dependence},
  journal = {Stochastics},
  volume = {94},
  number = {1},
  pages = {111--142},
  year = {2022},
  doi = {10.1080/17442508.2021.1914620}
}

@article{Marinucciv2,
  author = {Marinucci, D. and Vadlamani, S.},
  title = {High-frequency asymptotics for {Lipschitz--Killing} curvatures of excursion sets on the sphere},
  journal = {The Annals of Applied Probability},
  volume = {26},
  number = {1},
  pages = {462--506},
  year = {2016},
  doi = {10.1214/15-aap1097}
}

@article{MarinucciRV,
  author = {Marinucci, D. and Rossi, M. and Vidotto, A.},
  title = {Non-universal fluctuations of the empirical measure for isotropic stationary fields on {$\mathbb S^2\times\mathbb R$}},
  journal = {The Annals of Applied Probability},
  volume = {31},
  number = {5},
  pages = {2311--2349},
  year = {2021},
  doi = {10.1214/20-aap1648}
}

@book{NourdinPeccati2012,
  author = {Nourdin, I. and Peccati, G.},
  title = {Normal Approximations with {Malliavin} Calculus: From {Stein}'s Method to Universality},
  series = {Cambridge Tracts in Mathematics},
  volume = {192},
  year = {2012},
  publisher = {Cambridge University Press},
  address = {Cambridge},
  doi = {10.1017/cbo9781139084659}
}

@article{NourdinNualart2020,
  author = {Nourdin, I. and Nualart, D.},
  title = {The functional {Breuer--Major} theorem},
  journal = {Probability Theory and Related Fields},
  volume = {176},
  number = {1--2},
  pages = {203--218},
  year = {2020},
  doi = {10.1007/s00440-019-00917-1}
}

@article{OvalleRM24,
  author = {Ovalle-Mu{\~n}oz, D. P. and Ruiz-Medina, M. D.},
  title = {{LRD} spectral analysis of multifractional functional time series on manifolds},
  journal = {TEST},
  volume = {33},
  number = {2},
  pages = {564--588},
  year = {2024},
  doi = {10.1007/s11749-023-00913-7}
}

@book{PipirasTaqqu17,
  author = {Pipiras, V. and Taqqu, M. S.},
  title = {Long-Range Dependence and Self-Similarity},
  series = {Cambridge Series in Statistical and Probabilistic Mathematics},
  volume = {45},
  year = {2017},
  publisher = {Cambridge University Press},
  address = {Cambridge},
  doi = {10.1017/cbo9781139600347}
}

@article{RuizMedina2022,
  author = {Ruiz-Medina, M. D.},
  title = {Spectral analysis of multifractional {LRD} functional time series},
  journal = {Fractional Calculus and Applied Analysis},
  volume = {25},
  number = {4},
  pages = {1426--1458},
  year = {2022},
  doi = {10.1007/s13540-022-00053-z}
}

@article{Taqqu,
  author = {Taqqu, M. S.},
  title = {Weak convergence to fractional {Brownian} motion and to the {Rosenblatt} process},
  journal = {Zeitschrift f{\"u}r Wahrscheinlichkeitstheorie und Verwandte Gebiete},
  volume = {31},
  number = {4},
  pages = {287--302},
  year = {1975},
  doi = {10.1007/bf00532868}
}

@article{Taqqu2,
  author = {Taqqu, M. S.},
  title = {Convergence of integrated processes of arbitrary {Hermite} rank},
  journal = {Zeitschrift f{\"u}r Wahrscheinlichkeitstheorie und Verwandte Gebiete},
  volume = {50},
  number = {1},
  pages = {53--83},
  year = {1979},
  doi = {10.1007/bf00535674}
}

@article{GiraitisSurgailis1985,
  author = {Giraitis, Liudas and Surgailis, Donatas},
  title = {{CLT} and other limit theorems for functionals of {G}aussian processes},
  journal = {Zeitschrift f{\"u}r Wahrscheinlichkeitstheorie und Verwandte Gebiete},
  volume = {70},
  number = {2},
  pages = {191--212},
  year = {1985},
  doi = {10.1007/BF02451428}
}

@article{BaiTaqqu,
  author = {Bai, S. and Taqqu, M. S.},
  title = {Multivariate limit theorems in the context of long-range dependence},
  journal = {Journal of Time Series Analysis},
  volume = {34},
  number = {6},
  pages = {717--743},
  year = {2013},
  doi = {10.1111/jtsa.12046}
}

@article{Major2026Noncentral,
  author = {Major, P{\'e}ter},
  title = {Non-central limit theorem for non-linear functionals of vector valued {G}aussian stationary random fields},
  journal = {Stochastic Processes and their Applications},
  volume = {200},
  pages = {105005},
  year = {2026},
  publisher = {Elsevier},
  doi = {10.1016/j.spa.2026.105005},
  url = {https://www.sciencedirect.com/science/article/abs/pii/S0304414926001377}
}

@article{SanchezDeNaranjo1993,
  author = {S{\'a}nchez de Naranjo, M. V.},
  title = {Non-central limit theorems for non-linear functionals of {$k$} {G}aussian fields},
  journal = {Journal of Multivariate Analysis},
  volume = {44},
  number = {2},
  pages = {227--255},
  year = {1993},
  doi = {10.1006/jmva.1993.1013}
}

@article{Arcones1994,
  author = {Arcones, Miguel A.},
  title = {Limit theorems for nonlinear functionals of a stationary {G}aussian sequence of vectors},
  journal = {The Annals of Probability},
  volume = {22},
  number = {4},
  pages = {2242--2274},
  year = {1994},
  doi = {10.1214/aop/1176988503}
}

@misc{BetkenDehling2021,
  author = {Betken, Annika and Dehling, Herold},
  title = {Test for independence of long-range dependent time series using distance covariance},
  year = {2021},
  howpublished = {\href{https://arxiv.org/abs/2107.03041}{arXiv:2107.03041}},
  eprint = {2107.03041},
  archivePrefix = {arXiv},
  primaryClass = {math.ST},
  doi = {10.48550/arXiv.2107.03041},
  url = {https://arxiv.org/abs/2107.03041},
  note = {First posted in 2021; revised 27 January 2026 (v3)},
  archiveprefix = {arXiv},
  primaryclass = {math.ST}
}

@book{Bosq2000,
  author = {Bosq, Denis},
  title = {Linear Processes in Function Spaces: Theory and Applications},
  series = {Lecture Notes in Statistics},
  volume = {149},
  year = {2000},
  publisher = {Springer},
  address = {New York},
  doi = {10.1007/978-1-4612-1154-9}
}

@article{Merlevede1996,
  author = {Merlev{\`e}de, Florence},
  title = {Central limit theorem for linear processes with values in a {H}ilbert space},
  journal = {Stochastic Processes and their Applications},
  volume = {65},
  number = {1},
  pages = {103--114},
  year = {1996},
  doi = {10.1016/S0304-4149(96)00099-3}
}

@article{MerlevedePeligradUtev1997,
  author = {Merlev{\`e}de, Florence and Peligrad, Magda and Utev, Sergey},
  title = {Sharp conditions for the {CLT} of linear processes in a {H}ilbert space},
  journal = {Journal of Theoretical Probability},
  volume = {10},
  number = {3},
  pages = {681--693},
  year = {1997},
  doi = {10.1023/A:1022653728014}
}

@article{CharaciejusRackauskas2014,
  author = {Characiejus, Vaidotas and Ra{\v{c}}kauskas, Alfredas},
  title = {Operator self-similar processes and functional central limit theorems},
  journal = {Stochastic Processes and their Applications},
  volume = {124},
  number = {8},
  pages = {2605--2627},
  year = {2014},
  doi = {10.1016/j.spa.2014.03.007}
}

@article{Duker2018,
  author = {D{\"u}ker, Marie-Christine},
  title = {Limit theorems for {H}ilbert space-valued linear processes under long range dependence},
  journal = {Stochastic Processes and their Applications},
  volume = {128},
  number = {5},
  pages = {1439--1465},
  year = {2018},
  doi = {10.1016/j.spa.2017.07.015}
}

@article{LiRobinsonShang2020,
  author = {Li, Degui and Robinson, Peter M. and Shang, Han Lin},
  title = {Long-range dependent curve time series},
  journal = {Journal of the American Statistical Association},
  volume = {115},
  number = {530},
  pages = {957--971},
  year = {2020},
  doi = {10.1080/01621459.2019.1604362}
}

@misc{DukerZoubouloglou2025,
  author = {D{\"u}ker, Marie-Christine and Zoubouloglou, Pavlos},
  title = {Scaling limits for sample autocovariance operators of {H}ilbert space-valued linear processes},
  year = {2025},
  howpublished = {\href{https://arxiv.org/abs/2506.17168v1}{arXiv:2506.17168v1}},
  eprint = {2506.17168v1},
  archivePrefix = {arXiv},
  primaryClass = {math.PR},
  doi = {10.48550/arXiv.2506.17168},
  url = {https://arxiv.org/abs/2506.17168v1},
  archiveprefix = {arXiv},
  primaryclass = {math.PR}
}

@article{PanaretosTavakoli2013,
  author = {Panaretos, Victor M. and Tavakoli, Shahin},
  title = {Fourier analysis of stationary time series in function space},
  journal = {The Annals of Statistics},
  volume = {41},
  number = {2},
  pages = {568--603},
  year = {2013},
  doi = {10.1214/13-AOS1086}
}

\end{document}